\documentclass[3p, number,sort&compress,times,fleqn, review]{elsarticle}

\usepackage{amsmath}
\usepackage{amssymb}
\usepackage{mathtools}
\usepackage[colorlinks = true,
            linkcolor = blue,
            urlcolor  = blue,
            citecolor = blue,
            anchorcolor = blue]{hyperref}

\usepackage[utf8]{inputenc}
\usepackage{subfig}
\usepackage{graphicx}

\usepackage{footnote}

\usepackage[noend]{algorithmic}
\usepackage{algorithm}

\usepackage{./csml}

\graphicspath{{./figs/}}

\begin{document}

\begin{frontmatter}

\title{Isogeometric analysis using manifold-based smooth basis functions}

\author{M. Majeed}  
\author{F. Cirak\corref{cor1}}  
\ead{f.cirak@eng.cam.ac.uk}

\cortext[cor1]{Corresponding author}

\address{Department of Engineering, University of Cambridge, Trumpington Street, Cambridge CB2 1PZ, U.K.}

\begin{abstract}
We present an isogeometric analysis technique that builds on manifold-based
smooth basis functions for geometric modelling and analysis.  Manifold-based surface construction techniques are well known in geometric modelling and a
number of variants exist. Common to all is the concept of constructing a smooth
surface by blending together overlapping patches (or, charts), as in 
differential geometry description of manifolds. Each patch on the surface has
a corresponding planar patch with a smooth one-to-one mapping onto
the surface. In our implementation, manifold techniques are combined with
conformal parameterisations and the partition-of-unity method for  deriving
smooth basis functions on unstructured quadrilateral meshes.  Each vertex and
its adjacent elements on the surface control mesh have a corresponding planar
patch of elements. The star-shaped planar patch with congruent wedge-shaped
elements is smoothly parameterised with copies of a conformally mapped unit
square. The conformal maps can be easily inverted in order to compute the
transition functions between the different planar patches that have an overlap
on the surface.  On the collection of star-shaped planar patches the partition
of unity method is used for approximation. The smooth partition of unity, or
blending functions, are assembled from tensor-product b-spline segments defined
on a unit square.  On each patch a polynomial with a prescribed degree  is
used as a local approximant.  In order to obtain a mesh-based approximation scheme the coefficients of the  local approximants are expressed in dependence of vertex coefficients.  This yields a basis function for each vertex of the mesh which is smooth and non-zero over a vertex and its adjacent elements. Our numerical simulations indicate the optimal convergence of the resulting approximation scheme for Poisson problems and near optimal convergence for thin-plate and thin-shell problems discretised with structured and unstructured quadrilateral meshes.
\end{abstract}


\end{frontmatter}

\newpage
\excludecomment{editing}


\begin{editing}
	\textcolor{red}{
		\begin{itemize} 
			\item basis functions in terms of $eta$ or $xi$
		\end{itemize}
	}
\end{editing}

%
\section{Introduction}  

The interoperability limitation of Computer Aided Design (CAD) and Finite
Element Analysis (FEA) systems has become one of the  major bottlenecks in
simulation-based design. CAD and FEA are inherently incompatible because they use for
historical reasons different mathematical representations.  As
advocated in isogeometric analysis the use of identical basis functions for CAD
and FEA can facilitate their integration. Today most of the research on
isogeometric analysis focuses on NURBS~\cite{Hughes:2005aa,Cottrell:2009aa} and
the related t-splines~\cite{Bazilevs:2010aa} and subdivision basis
functions~\cite{Cirak:2002aa}. The inherent tensor-product structure of NURBS
means that additional techniques are required for  geometries that are composed
out of several NURBS patches.  Specifically, around extraordinary (or irregular)
points where the number of patches that join together is different than four,
i.e. $ v \neq 4$, alternative techniques are necessary to maintain smoothness.
One prevalent approach in geometric design is to introduce additional higher
order patches around the extraordinary point and to ensure that all patches
match up $G^k$ continuously at their boundaries.   $G^k$ refers to the notion of
geometric continuity and, for instance, $G^1$ implying  tangent plane
continuity. As recently pointed out by Groisser et
al.~\cite{groisser2015matched} in isogeometric analysis $G^k$ leads to $C^k$
continuity because the geometry and field variables are interpolated with the
same basis functions. The utility of classical $G^k$ constructions in
isogeometric analysis has recently been demonstrated in a number of
papers~\cite{nguyen2014comparative,nguyen2016c}.  $G^k$ constructions have  also
been explored in the context of isogeometric analysis with
t-splines~\cite{scott2013isogeometric}. A different approach for dealing with
extraordinary points is provided by subdivision surfaces. The neighbourhood of
the extraordinary point is replaced by a sequence of nested $C^k$ continuous
patches which join $C^1$ continuously at the point
itself~\cite{Peters:2008aa,Zorin:2000aa}. Subdivision basis functions for finite
element analysis have originally been proposed  in~\cite{Cirak:2000aa} and  have
been more intensely studied  in a number of recent
papers~\cite{Cirak:2011aa,juttler2016numerical,wei2015truncated}.

\begin{figure}
\centering
\begin{minipage}{.27\textwidth}
	\centering
	\includegraphics[scale=0.4]{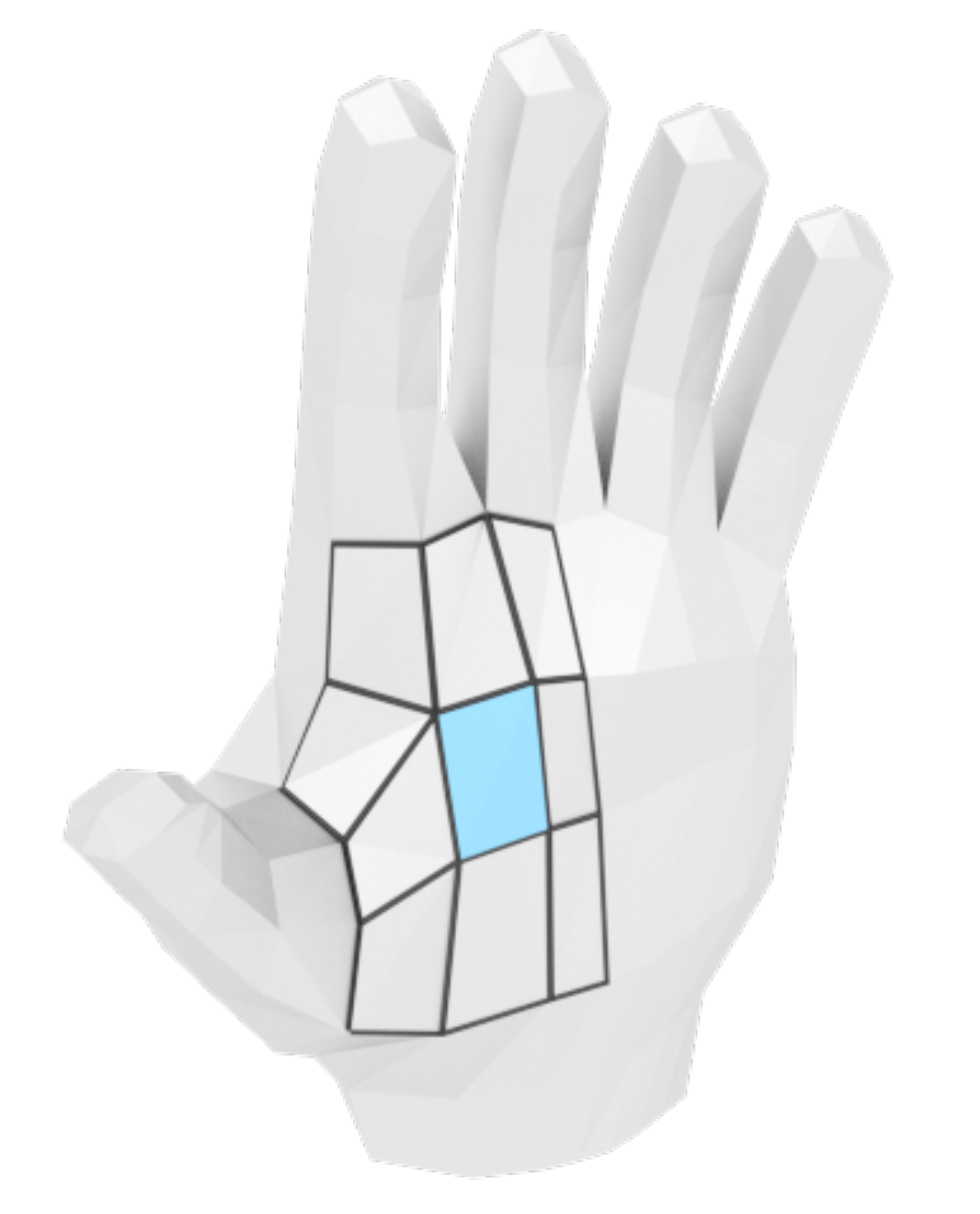}
\end{minipage} %
 \hspace*{\fill} 
\begin{minipage}{.27\textwidth}
	\centering
	\includegraphics[scale=0.4]{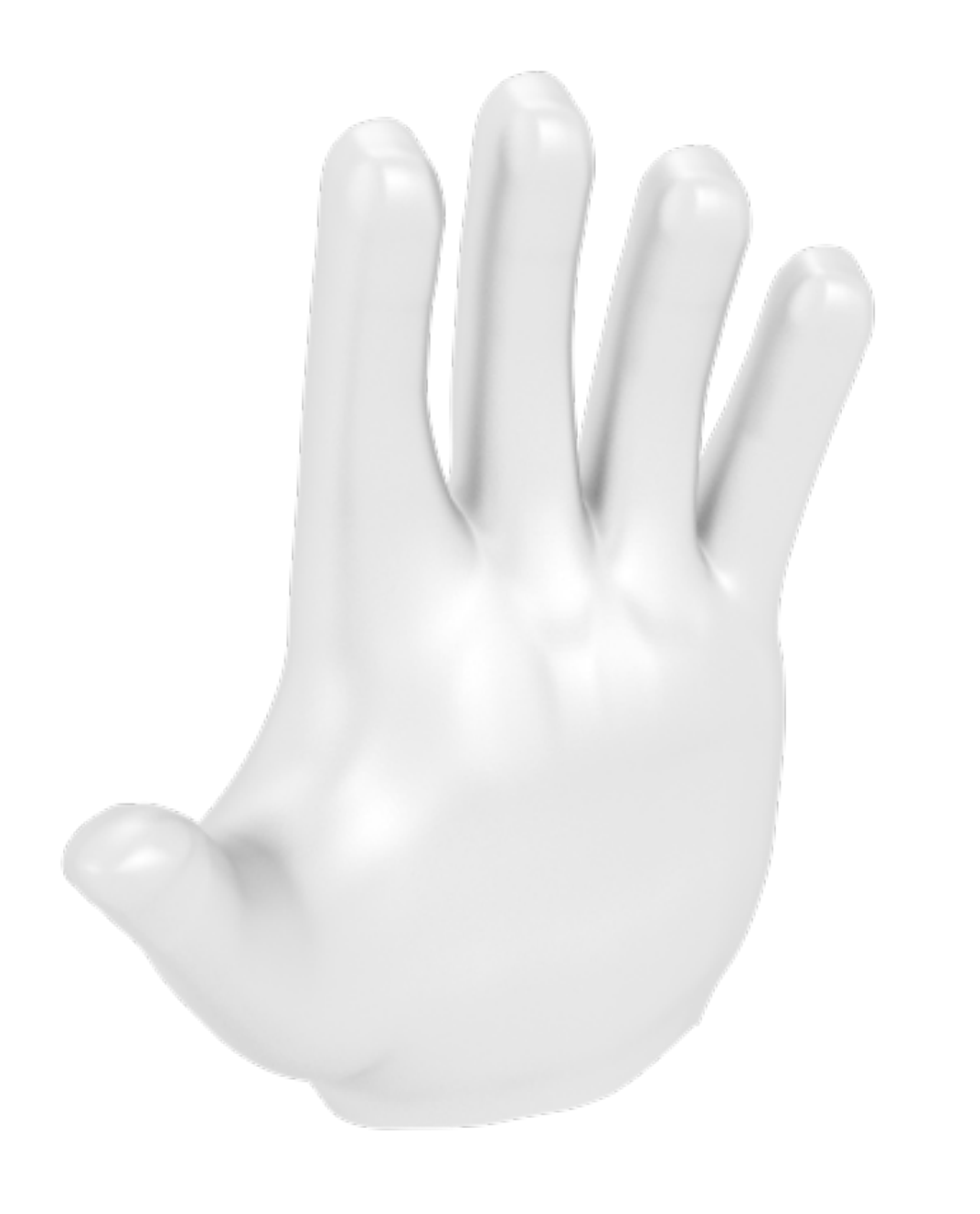}
\end{minipage}
 \hspace*{\fill} 
\begin{minipage}{.36\textwidth}
	\centering
	\includegraphics[scale=0.4]{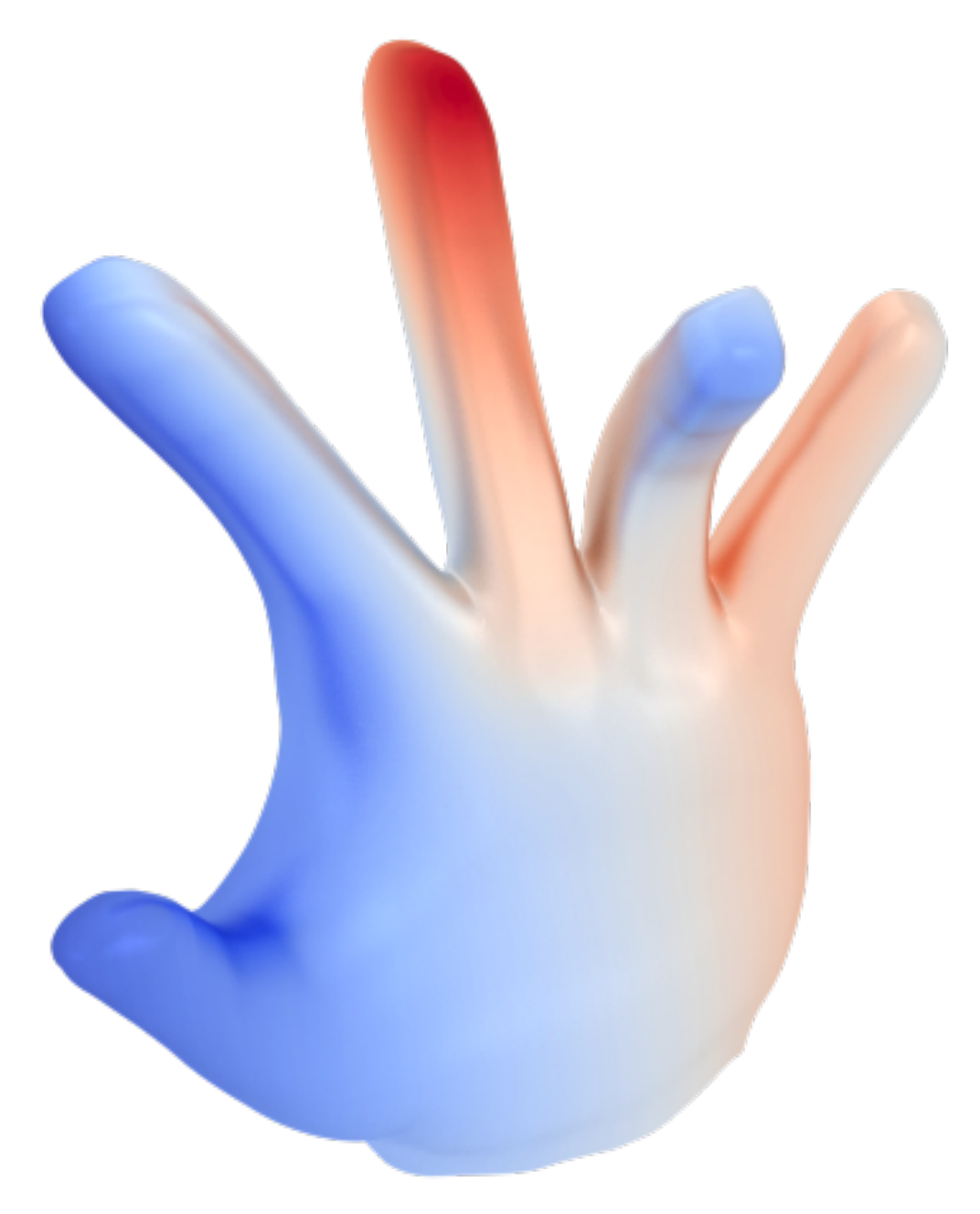}
\end{minipage}
\caption{Isogeometric analysis using manifold-based smooth surfaces. The smooth surface in the middle is described by the control mesh on the left and the corresponding manifold basis functions.  The deflected geometry on the right is the result of a finite element thin-shell computation with manifold basis functions.} \label{fig:introducEx}
\end{figure}

We introduce in this paper an isogeometric analysis  technique that builds on
manifold-based basis functions for geometric modelling and analysis. As known
from differential geometry, manifolds provide a rigorous framework for
describing and analysing surfaces with arbitrary topology;
see~\cite{doCarmo1976differential, schutz1980geometrical}. Informally, with
manifolds a surface in Euclidean space $\varmathbb{R}^3$ is constructed from
mapping and blending together planar patches from $\varmathbb{R}^2$.  Manifold
techniques for mesh-based construction of  smooth $C^k$ continuous surfaces were
first considered in Grimm et al.~\cite{grimm1995modeling}.   Other mesh-based
manifold constructions have later been proposed, e.g.,
in~\cite{navau2000modeling,ying2004simple,gu2006manifold,della2008construction}.
In a manner similar to the description of splines, a $C^k$ continuous surface is
described with a quadrilateral or triangular control mesh and each vertex has a corresponding basis function
with a local support, see Figure~\ref{fig:introducEx}.
In contrast to the aforementioned $G^k$ constructions which rely on matching up separate patches,
in manifolds a $C^k$ continuous surface is created by smoothly blending  of
overlapping patches. The idea of blending surfaces from overlapping patches is a
common theme in geometric modelling and has been used, for instance,  for
increasing the smoothness of subdivision surfaces around the extraordinary
vertices~\cite{levin2006modified, zorin2006constructing,
antonelli2013subdivision} or (meshfree) point-based surface
processing~\cite{levin2004mesh, pauly2002efficient}. In Millan et
al.~\cite{Millan:2010aa,millan2013nonlinear}  point-based surface blending
techniques have also been used for meshfree thin-shell analysis.

In the present work we follow Ying and Zorin~\cite{ying2004simple} and construct smooth basis functions by combining manifold techniques with conformal parameterisations and partition of unity method. The control mesh consists of quadrilateral elements with some extraordinary vertices (i.e. $ v \neq 4$ for some non-boundary vertices) and the construction gives one basis function for each vertex. 
The first step is to assign  each vertex of the control mesh and its adjacent elements  a planar sub-mesh with the same connectivity. The sub-meshes serve as control meshes for planar surface patches, which can be understood as parameter spaces for basis functions. For $C^k$ continuous basis functions  the planar patches  have to have a $C^k$ smooth parameterisation.  Although other choices are conceivable, the patches are parameterised using  conformal (angle-preserving) maps.  Since each surface point is represented on several patches transition functions composed of conformal maps are used  to navigate between adjacent patches. In the second step of the construction,  on each  planar patch  the conventional partition of unity method (PUM) of Melenk et al.~\cite{melenk1995generalized,Melenk1996a}  is used for constructing  basis functions.  According to PUM, the basis functions  are the product of a partition of unity function and a patch specific polynomial approximant. In computer graphics literature the partition of unity function and the patch specific polynomial basis are usually referred to as the blending function and the embedding function, respectively. We use as partition of unity functions b-splines that have zero value and zero $k-1$ derivatives at the patch boundaries. In order to enforce partition of unity the b-splines on different patches overlapping the same point on the manifold surface are first identified with transition functions and subsequently normalised as in usual PUM. The last step in the basis function  construction is to express the local PUM polynomial approximant in dependence of vertex values using a least-squares approximation.  In this mesh-based approach the degree of the polynomial approximant and the number of vertices in a sub-mesh are correlated. In order to increase the polynomial degree the sub-meshes are enlarged with mesh refinement by quadrisectioning.  The basis functions depend only on the connectivity of the control mesh but not its geometry so that they can be  precomputed and tabulated for different valences $v$. 

The outline of this paper is as follows. Section 2  reviews the relevant manifold concepts from differential geometry and the partition of unity method. The mesh-based manifold basis functions are introduced in Section 3. First, one-dimensional polygonal control meshes are considered, even though  
 it is straightforward to combine one-dimensional manifolds with the partition of unity method by using simple transition functions. Subsequently, two-dimensional quadrilateral meshes are considered  for which conformal maps are used as transition functions. In both dimensions, it is shown how the polynomial coefficients in the partition of unity approximation can be expressed as a linear combination of vertex coefficients. In Section 4 the derived mesh-based basis functions are applied to a number of Poisson and thin-shell problems. Numerical convergence with increasing mesh size on meshes with and without extraordinary vertices is investigated.

	%
\section{Preliminaries}
%

%
\subsection{Review of manifold concepts \label{sec:manifold}}
%
In the following we provide an informal introduction to differentiable manifolds
with the aim to introduce the necessary terminology. For clarity, our discussion
is restricted to surfaces, i.e. two-dimensional manifolds, embedded in the
three-dimensional Euclidean space $\varmathbb{R}^3$. A similar introduction, but
oriented more towards geometric modelling, can be found in~\cite{grimm2006surface}. The
manifold concept is much more general than needed in this paper. A more rigorous
discussion is found in most standard differential geometry textbooks, see,
e.g.,~\cite{doCarmo1976differential,schutz1980geometrical}.

A regular  two-manifold, or surface, $\Omega$ is defined as the set of points in  $\varmathbb{R}^3$, which can be locally continuously one-to-one mapped onto a set of points  in $\varmathbb{R}^2$. This definition naturally  excludes surfaces with t-joints or isolated points. In applications this is, however, not a major restriction since a geometry, for instance, with a t-joint can be represented as a collection of several manifolds.

\begin{figure}
  \centering
  \includegraphics[scale=1.2]{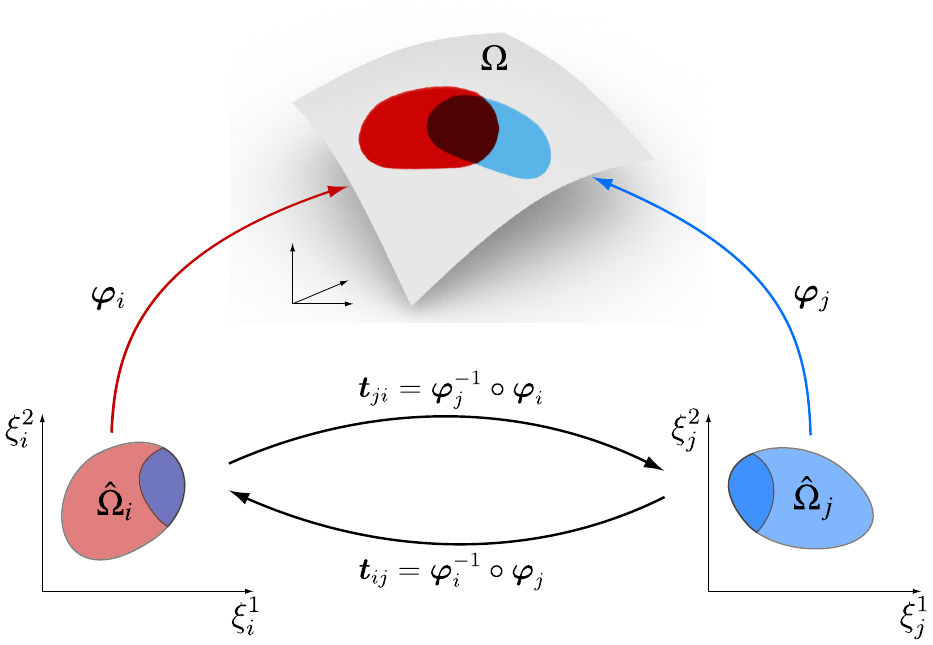}
  \caption{A manifold $\Omega$ and two of the charts $(\hat{\Omega}_i, \vec{\varphi}_i)$ and $(\hat{\Omega}_j, \vec{\varphi}_j)$ for representing it. A point $P$ in the dark shaded area of the manifold $\Omega$ is represented in both charts and has different coordinates in $\hat{\Omega}_i$ and $\hat{\Omega}_j$. The change of coordinates between the two charts is given by the  transition maps $\vec t_{ij} = \vec t_{ji}^{-1}$. \label{fig:mc_1}}
\end{figure}

By definition on a regular surface $\Omega$ around each point $P$ there is an open region  $\Omega_i \in \Omega$ that can be mapped to an open planar region $\hat{\Omega}_i \in \varmathbb{R}^2$, see Figure~\ref{fig:mc_1}. In line with the  partition of unity method  terminology, in the following ${\Omega}_i$ will be referred to as a patch and $\hat{\Omega}_i$ as a planar patch. Furthermore, we denote the function for mapping between two patches with  $\vec{\varphi}_i : \hat{\Omega}_i \rightarrow \Omega_i $. On the surface $\Omega$ there are many overlapping patches $\Omega_i$ such that
\begin{align}
 		 \Omega &= \bigcup_{i \in \varmathbb{N}} \Omega_i 
\end{align}
and each point $P$  lies at least on one patch. The pair consisting of $(\hat{\Omega}_i, \vec{\varphi}_i)$ is called a chart. The set of all charts $\{(\hat{\Omega}_i, \vec{\varphi}_i) \}$ is referred to as an atlas for representing the surface $ \Omega$. 
As illustrated in Figure~\ref{fig:mc_1}, each planar patch $\hat{\Omega}_i$ has its own coordinate system $\vec {\xi}_i = (\xi^1_i, \xi^2_i)$. The same point $P \in \Omega$ in  the intersection between the two patches $ P \in ( \Omega_i \cap \Omega_j)$ has the coordinates $\vec {\xi}_i = (\xi^1_i, \xi^2_i)$ in $\hat{\Omega}_i$  and $\vec {\xi}_j = (\xi^1_j, \xi^2_j)$ in $\hat{\Omega}_j$.  In order to compute the underlying coordinate transformations  we introduce the  transition functions
\begin{equation}
	\vec t_{ji} = \vec{\varphi}_j^{-1} \circ \vec {\varphi}_i
\end{equation}
that are composed out of the mappings $\vec{\varphi}_j$. The transition functions are symmetric $\vec t_{ij} = \vec t_{ji}$ and satisfy the cocyle condition $\vec t_{ij} \circ \vec t_{ki} = \vec t_{kj}$ when the preimage of $P$ lies in three planar patches $\hat \Omega_i$, $\hat \Omega_j$ and $\hat \Omega_k$. For a surface to be $C^k$ continuous the transition maps must be $C^k$ continuous. Evidently, for a differentiable surface $k$ has to be equal or larger than one.

%
\subsection{Review of the Partition of Unity Method (PUM)}
The  first use of PUM for creating finite element basis functions goes back to
the seminal work of  Babuska et al.~\cite{babuvska1994special} and was
subsequently further developed, for instance, in Melenk et
al.~\cite{Melenk1996a} and Duarte et al.~\cite{duarte1996h}.  In PUM a given
domain $\Omega$ in the Euclidean space \mbox{$\varmathbb{R}^d$}, with $ d\in
\{1,2,3 \}$, is partitioned into overlapping patches $\Omega_i$ such that
\begin{align}
 		 \Omega &= \bigcup_{i \in \varmathbb{N}} \Omega_i  \, .
\end{align}
As opposed to the manifolds introduced in previous Section~\ref{sec:manifold}, there is only one single coordinate system in the Euclidean space~$\varmathbb{R}^d$ and a point~$P \in \Omega$ has the same coordinates on~$\Omega$ and all patches~$\Omega_i$. Hence, the transition functions between the different patches are identity maps.

Next, a blending function $w_i$ is defined on each patch $\Omega_i$.  By definition the sum of the blending functions $w_i$ over all the regions is
\begin{equation}
 	\sum_i w_i  = 1 \, .
\end{equation}
The set of blending functions~$\{ w_i \}$ is also referred to as the partition of unity subordinate to the set of open patches~$\{ \Omega_i \}$.  As will become clear, in order to obtain smooth PUM basis functions the blending functions $w_i$ have to be smooth.  In addition, on the region boundaries~$\partial \Omega_i$  the function value and derivatives of $w_i$  have to be zero. After choosing on each patch an arbitrary functions~$\hat w_i$ that have the prerequisite properties they can be normalised to yield a blending function 
\begin{equation} \label{eq:puScaling}
		w_i = \frac{\hat{w}_i }{\sum_j \hat{w}_j }  \, .
\end{equation}
In our applications $\hat w_i$ are usually b-spline basis functions. 

On each patch~$\Omega_i$, in addition to the blending functions $w_i$ a local
polynomial approximant $f_i$ is considered
\begin{equation}
	f_i = \vec \alpha_i \cdot \vec p_i \,  ,
\end{equation}
where~${\vec p}_i$ is vector containing a complete polynomial basis, $\vec {\alpha}_i$ is the corresponding vector of  the coefficients and the dot represents their scalar product. For instance, for one-dimensional domains and a monomial basis the two vectors are of the form
\begin{editing}
\todo{shouldn't $\left (x^1 \right)^1 \ldots$ be replaced by $\left
(\xi^1 \right)^1$ in equation 8}
\end{editing}
\begin{align}
	\vec \alpha &= \begin{bmatrix} \alpha_1 & \alpha_2 & \alpha_3 & \alpha_4 &  \cdots \end{bmatrix} ^\trans  \, \\
	\vec p  &= \begin{bmatrix} 1& \left (x^1 \right)^1 & \left (x^1 \right)^2  & \left (x^1 \right)^3 & \cdots  \end{bmatrix} ^\trans  \, .
\end{align}

The global approximant $f$ is the sum of the local approximants~$f_i$  and their
multiplication with the blending functions~$w_i$, that is,
\begin{align}  \label{eq:pumSummation}
	f  = \sum_i w_i  f_i  = \sum_i w_i  \left ( \vec   \alpha_i \cdot \vec p_i \right )  \, . 
\end{align}
The smoothness of this function depends on the smoothness of $w_i$ and $f_i$. Since the polynomials in $f_i$ are infinitely smooth,  the smoothness of $f$  is exclusively controlled by the blending functions $w_i$. According to the convergence estimates given in~\cite{Melenk1996a}, the convergence  rates for the approximant $f$ depend on the degree of the polynomial basis $\vec {p}_i$ and the constants on the blending functions $ w_i$ and the layout of the overlaps $\Omega_i$.

\begin{figure}[]
  \centering
  	\begin{minipage}{.49\textwidth}
		\centering
		\includegraphics[scale=1]{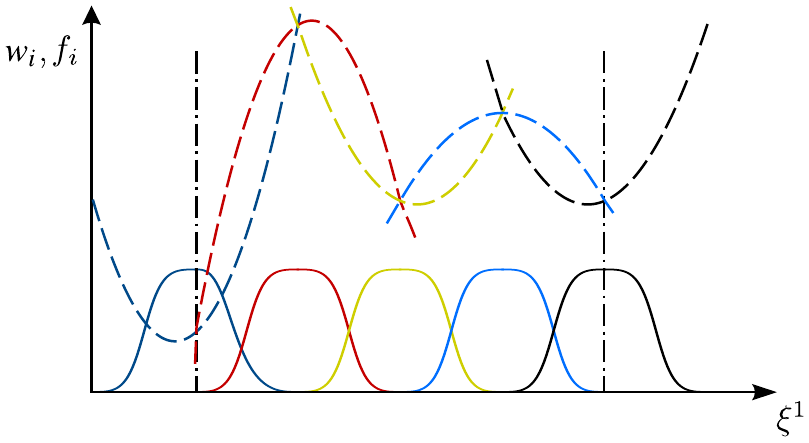}
	\end{minipage}
 	\hspace*{\fill} 
	\begin{minipage}{.49\textwidth}
		\centering
		\includegraphics[scale=1]{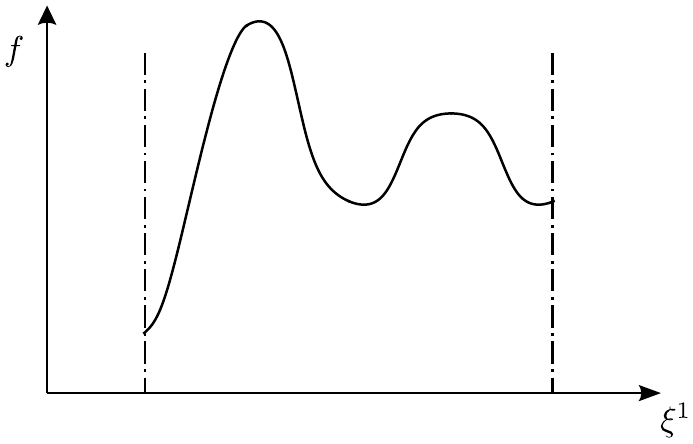}
	\end{minipage}
    \caption{Construction of a smooth function with PUM on a domain with five
  overlapping patches. On the left the blending functions $w_i$ (solid) and the local polynomials $f_i$ (dashed) are shown.  On the right the obtained smooth function $f$ is shown.  \label{fig:smoothBasisFunctions_1} }
\end{figure}

The illustrative one-dimensional example in
Figure~\ref{fig:smoothBasisFunctions_1} showcases the construction of a smooth
function using PUM. On each patch~$\Omega_i$ one single cubic b-spline basis
functions is used as~$\hat w_i$ for defining the (normalised) blending function~$w_i$
according to~\eqref{eq:puScaling}. The support size of a  cubic b-spline~$\hat
w_i$ is  equal to the size  of its corresponding patch~$\Omega_i$. The cubic
b-spline is a $C^2$ continuous function and at the boundaries of~$\Omega_i$  its
value and first and second derivatives are zero~\cite{Farin:2002aa}. In addition
to the blending functions~$w_i$, in Figure~\ref{fig:smoothBasisFunctions_1} the
local quadratic polynomials $f_i$ are also shown. The  smooth function $f$
constructed according to~\eqref{eq:pumSummation} is shown in
Figure~\ref{fig:smoothBasisFunctions_1}(b). It is evident that this curve is by
construction~$C^2$ continuous.

\begin{figure}[]
  \centering
  \includegraphics[scale=0.95]{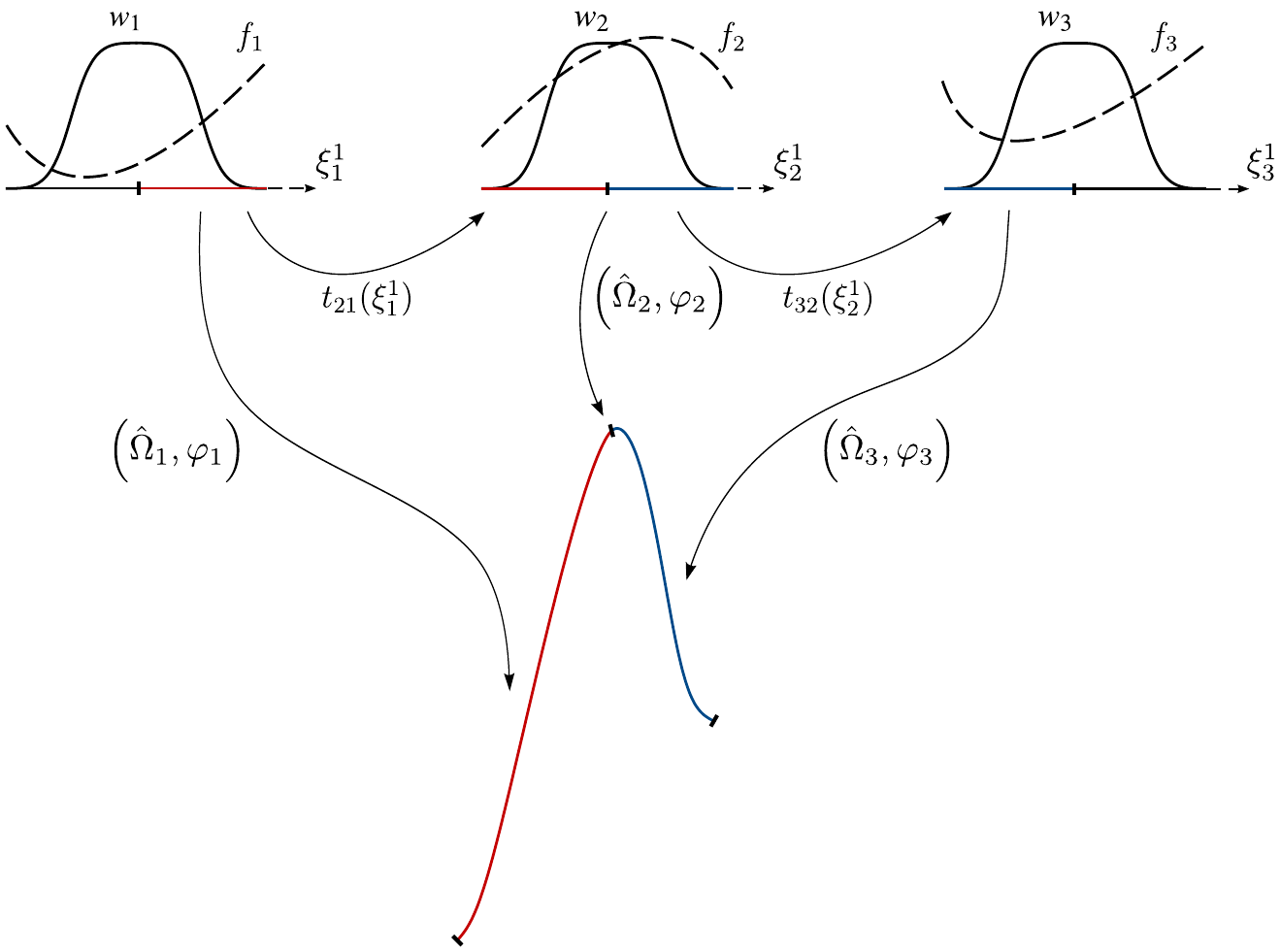}
  \caption{Construction of a smooth function with  manifold-based PUM. The
  domain consists of  three overlapping patches (top). On each patch the blending
  function $w_i$ (solid) and the local polynomial $f_i$ (dashed) are shown. The
  transition functions between the overlapping patches are chosen as~$t_{21}
  = - \left (\xi_1^1 \right )^2 + 2 \xi_1^1 - 1$ and~$t_{32} = \xi_2^1 - 1$.  The obtained smooth function $f$ is shown at the bottom.
  \label{fig:smoothBasisManifoldBased} }
\end{figure}

For the considered one-dimensional example, it is straightforward to use a manifold-based approach for constructing the partition of unity function. As shown in Figure~\ref{fig:smoothBasisManifoldBased}, the coordinate systems on each patch can be chosen differently. For switching between the different coordinate systems the transition functions~$t_{ij}$ are used. In this example one of the chosen transition functions, i.e. ~$t_{21}$, is quadratic and the other, i.e.~$t_{32}$ ,is linear. It is however possible to chose any other monotone~$C^2$ function. Choosing other transition functions will lead to a change in the shape of the constructed function~$f$.

%
\section{Mesh-based manifold basis functions} 
%
We are now in a position to introduce the construction of manifold-based basis functions on  one- and two-dimensional meshes. The idea of using manifolds for smooth interpolation on meshes has been originally introduced in computer graphics by Grimm et al.~\cite{grimm1995modeling}. The one-dimensional case is straightforward and is only discussed in order to provide some intuition for the two-dimensional case. The essential difficulty in two-dimensions lies in defining suitable transition functions.   We use the conformal maps  as introduced in  Ying and Zorin~\cite{ying2004simple} for defining the transition functions. Alternative definitions have been provided in ~\cite{grimm1995modeling,navau2000modeling, antonelli2013subdivision}. 

%
\subsection{One-dimensional meshes}
\begin{figure}
  \centering
  \includegraphics[scale=1]{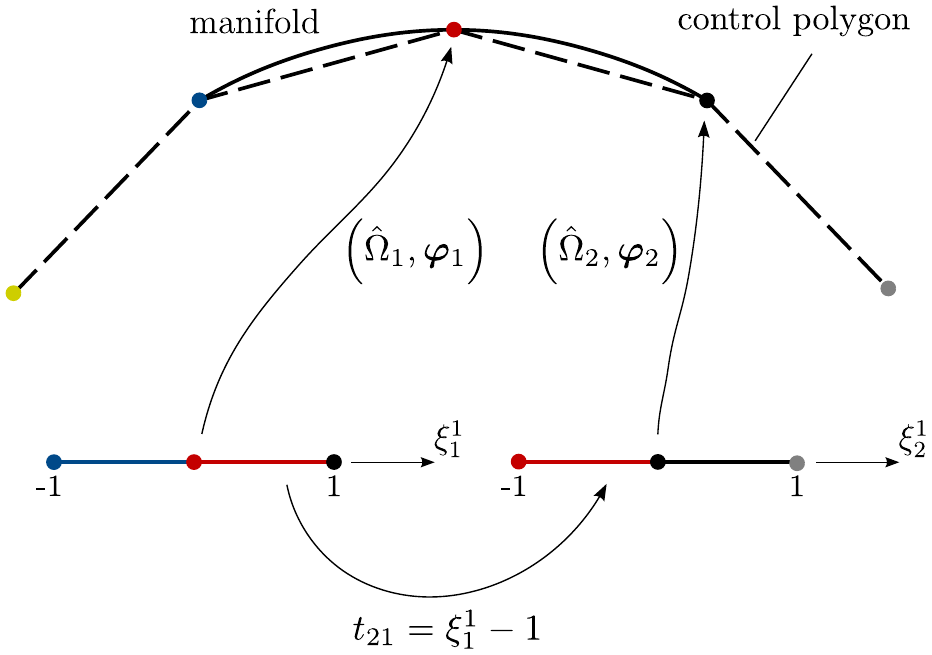}
  \caption{A control polygon (dashed) with the corresponding smooth curve (solid). On each of the two shown one-rings the blending function is a normalised cubic b-spline and the local approximant is a quadratic polynomial.}
  \label{fig:mbc_1}
\end{figure}

%
First we aim to construct a smooth curve, i.e. a one-dimensional smooth manifold for a given coarse control polygon in the Euclidean space~$\varmathbb{R}^3$.  We begin with defining charts~$(\hat{\Omega}_i, \vec\varphi_i)$ for each vertex, see~ Figure \ref{fig:mbc_1}. The planar patch~$\hat{\Omega}_i$ is formed from two segments and the attached three vertices. This is usually referred to as the one-ring of the centre vertex. It is also possible to increase the size of~$\hat{\Omega}_i$  to a two-ring or even larger.  The chosen size of the patches influences the number of overlapping patches at each point. As typical for manifolds, each planar patch $\hat{\Omega}_i$ has its own coordinate system. The scalar transition functions $t_{ij}$ enable to navigate  between the patches~$\hat{\Omega}_i$ and~$\hat{\Omega}_j$. The transition functions are chosen as linear maps.

The coordinates of a point  $P(x^1,x^2, x^3) $ on the smooth curve  is now
determined with the partition of unity method.  According
to~\eqref{eq:pumSummation} we can write, for instance, for the~$x^1$ coordinate
of the  point~$P$ with a preimage  on the planar patch~$\hat{\Omega}_i$ and
coordinate~$\xi_i$
\begin{equation} \label{eq:coordinatePUM}
	x^1 (\xi_i) = \sum_j w_j(\xi_i) \left (  \vec \alpha_j^1 \cdot \vec p (\xi_i)  \right ) \, .
\end{equation}
The summation is over all patches and in order to evaluate the sum it is
necessary to use the transition maps. Although we give here and in the following
only the expression for~$x^1$, the other two coordinates~$x^2$ and~$x^3$ are expressed similarly.  On each one-ring there are three vertices, which motivates the choice of a quadratic basis for~$\vec p (\xi_i)$. For a quadratic Lagrangian basis~$\vec L (\xi_i)$ the three coefficients are simply  the coordinates of the vertices in the one-ring:
\begin{equation} \label{eq:coordinatePUM0}
	x^1 (\xi_i) = \sum_j w_j(\xi_i) \left (  \vec x_j^1 \cdot \vec L  (\xi_i)  \right ) \, .
\end{equation}
Both vectors~$\vec x_j^1$ and~$\vec L (\xi_i)$ have three components and the entries of~$\vec x_j^1$ are the $x^1$ coordinates of the three control polygon  vertices. This equation is illustrated in Figure~\ref{fig:mbc_1}.  As can be seen, the smooth curve passes exactly through the vertices of the control polygon.   

\begin{figure}
\centering
	\begin{minipage}{.49\textwidth}
		\centering
		\includegraphics[scale=1.05]{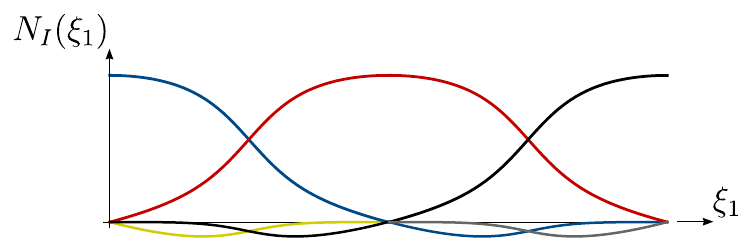}
	\end{minipage}
 	\hspace*{\fill} 
	\begin{minipage}{.49\textwidth}
		\centering
		\includegraphics[scale=1.05]{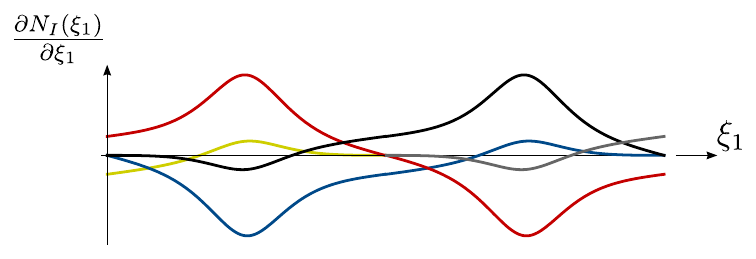}
	\end{minipage}
	\caption{ Basis functions (left) and their derivatives (right)  in
	two adjacent elements. In one element there are four  non-zero basis functions with
	two, i.e. the negative ones in the left figure, corresponding to vertices  on
	neighbouring elements.}	\label{fig:mbc_2}
\end{figure}

Next, we rewrite equation~\eqref{eq:coordinatePUM0}  in index notation to define
basis functions that can be used for finite element analysis:
\begin{equation} \label{eq:coordinatePUM1}
	x^1 (\xi_i)  = \sum_j w_j(\xi_i) \left (  \sum_{I=1}^3  x_{j,I}^1  L_I (\xi_i)  \right ) 
		              =   \sum_{I=1}^3  \left ( \sum_j w_j(\xi_i)  L_I (\xi_i)  \right ) x_{j,I}^1  =  \sum_{I=1}^3 N_I(\xi_i)  x_{j,I}^1 \, ,
\end{equation}
where $N_I$ are the three basis functions corresponding to the three vertices in the patch~$\hat{\Omega}_j$. In Figure~\ref{fig:mbc_2}   the non-zero basis functions and their derivatives in one patch are shown.  In the underlying construction the blending functions are normalised cubic b-splines, local polynomials are quadratic and the transition function are linear.  Note that the support size of one basis function is two elements. Due to the overlaps there are four non-zero basis functions in one  element. The resulting basis functions are $C^2$ continuous. 

%
\subsection{Two-dimensional quadrilateral meshes \label{sec:twoDconstruction}}
%
We now consider the construction of a smooth surface, i.e. smooth two-manifold, for a given coarse control mesh. Our approach follows  closely the construction originally introduced in Ying and Zorin~\cite{ying2004simple}. Although  only quadrilateral meshes are considered, it is straightforward to extend the technique to triangular meshes.

\begin{figure}
\centering
  	\begin{minipage}{.32\linewidth}    	
		\centering
		\includegraphics[scale=0.77]{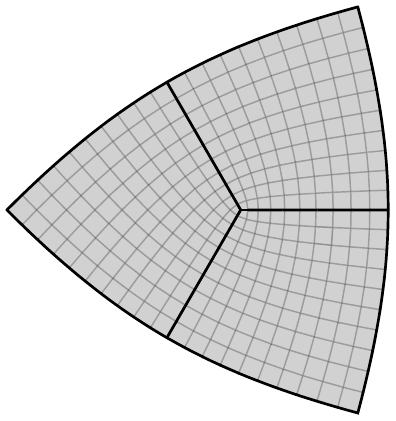}
	\end{minipage}
 	\hspace*{\fill} 
	\begin{minipage}{.32\linewidth}
		\centering
		\includegraphics[scale=0.77]{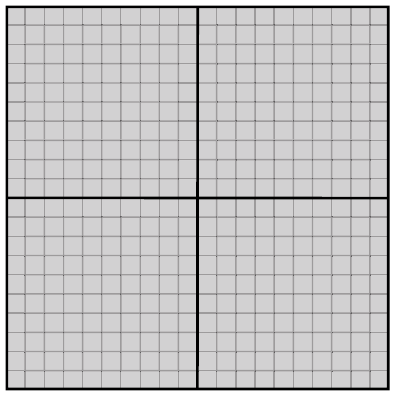}
	\end{minipage}
 	\hspace*{\fill} 
	\begin{minipage}{.32\linewidth}
		\centering
		\includegraphics[scale=0.77]{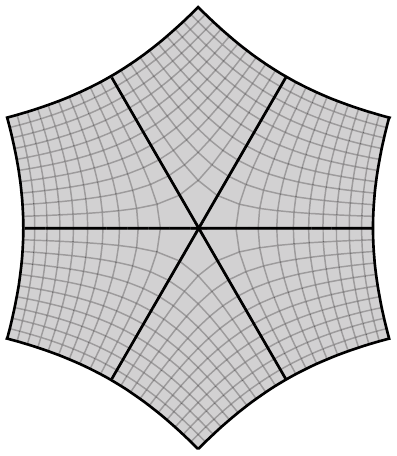}
	\end{minipage}
  \caption{Conformal parameterisations one one-ring patches with valence three (left), four  (middle) and six (right). The shown fine meshes  represent the parameter lines with either $\eta^1=\text{const.}$ or $\eta^2 = \text{const.}$.}
  \label{fig:mbc_4}
\end{figure}

\begin{figure}
  \centering
  \includegraphics[scale=0.95]{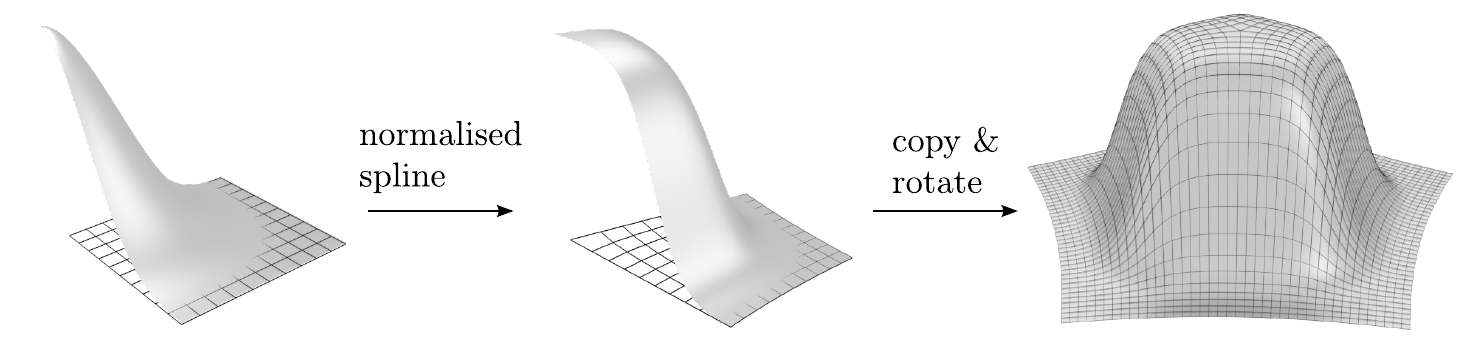}
  \caption{One quarter  of a tensor-product cubic b-spline as the non-normalised
   blending function $\hat w_i$ on the unit square (left). Normalised blending
   function $w_i$ on the conformal wedge (middle) and on a conformal patch with valence five (right).}
  \label{fig:mbc_5}
\end{figure}

\begin{figure}
\centering
	\begin{minipage}{.31\textwidth}
		\centering
		\includegraphics[scale=0.8]{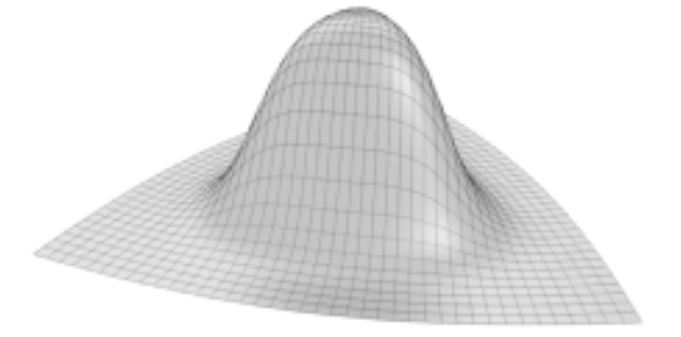}
	\end{minipage}
 	\hspace*{0.03\textwidth} 
	\begin{minipage}{.31\textwidth}
		\centering
		\includegraphics[scale=0.8]{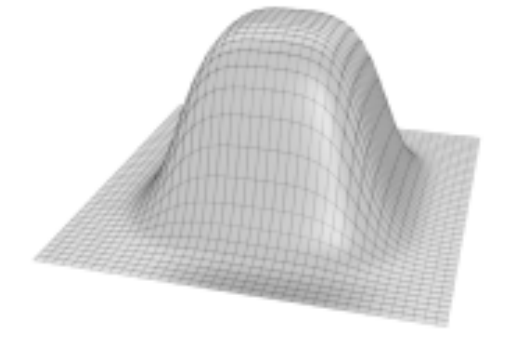}
		\end{minipage}
	\begin{minipage}{.31\textwidth}
		\centering
		\includegraphics[scale=0.8]{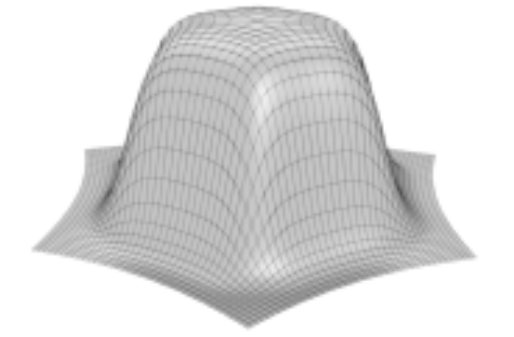}
		\end{minipage}
  \caption{Blending functions $w_i$ for one-ring patches with valence three (left), four (middle) and six (right). The shown fine meshes  represent the parameter lines with either $\eta^1=\text{const.}$ or $\eta^2 = \text{const.}$.}
  \label{fig:mbc_6}
\end{figure}

In addition to the smoothness properties of blending functions and local
polynomials, the smoothness of transition functions is central in generating
smooth surfaces. For implementation purposes, it is also important that the
transition functions and their inverses are readily computable. Similar to the
one-dimensional construction, we define charts~$(\hat{\Omega}_i,
\vec{\varphi_i})$  for each vertex of the mesh. The planar patch~$\hat{\Omega}_i$ is chosen for now as the one-ring of elements around a vertex. The
number of elements in the one-ring of a vertex is referred to as the valence~$v$
of the vertex. On structured meshes all vertices inside the domain have  valence~$v=4$ and on unstructured meshes it can be arbitrary. Hence, in the unstructured
case the overlapping one-ring patches can have different valences, which makes
the computation of a smooth transition function challenging.
Ying et al.~\cite{ying2004simple} proposed conformal maps as smooth and easy
computable transition functions. Recall here that conformal mapping is an angle
preserving transformation. The generation of a conformal parameterisation for a
one-ring of elements proceeds in several steps. In Figure~\ref{fig:mbc_3} the
procedure for a vertex  with valence~$v=5$ is illustrated. The smooth
parameterisation is obtained by conformally mapping, rotating and combining unit
squares. The points of the unit square have the coordinates~$\vec \eta =
(\eta^1, \eta^2)$ and are expressed as a complex number~$z=\eta^1+i \eta^2$. The
conformal transformation~$z^{4/v}$  maps the square to a wedge. In computing the
mapping recall the following standard relations:
\begin{figure}[]
  \centering
  	\includegraphics[scale=1.0]{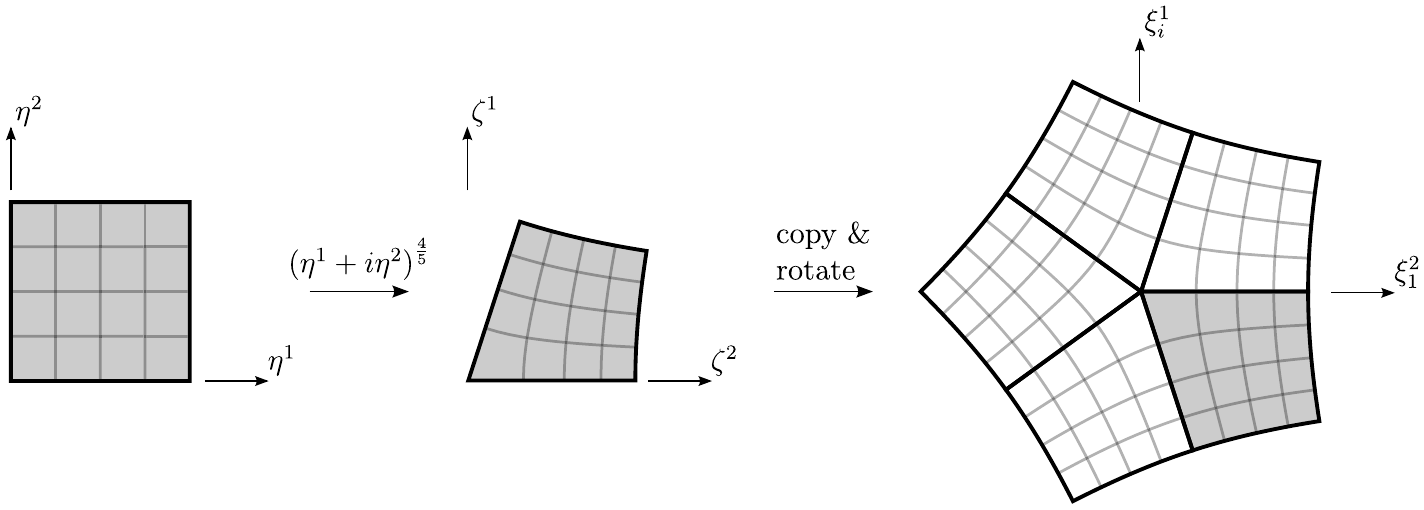}
  \caption{Conformal parameterisation of an one-ring patch with valence $v=5$ (right). The unit square (left) is first conformally mapped to a wedge (middle) with $z^{4/5} = (\eta^1 + i \eta^2)^{4/5}$. Subsequently five copies of the wedge are suitably rotated and combined into a one-ring patch.   \label{fig:mbc_3}}
\end{figure}
\begin{equation}
	 z = \eta^1 + i \eta^2 = |z|  \left ( \cos \phi + i \sin \phi  \right ) = | z | e^{i \phi}  
 \qquad \text{with } \; 
	|z| = \sqrt{{(\eta^1)}^2 + {(\eta^2)}^2} \; \text{ and } \; \phi = \arctan \left ( \eta^2 / \eta^1 \right )  \, .
\end{equation}
After the mapping the coordinates $\vec \eta = (\eta^1, \eta^2)$  of a point
become $\vec \zeta = (\zeta^1, \zeta^2)$ according to:
\begin{equation}
	 \zeta^1 + i \zeta^2 = z^{4/v} = |z|^{4/v}  e^{i 4\phi /v} \, .
\end{equation}
The wedge-shaped image of the conformally mapped unit square forms one sector of the one-ring patch. This wedge is copied and suitably rotated to form a smoothly parameterised one-ring forming patch~$\hat{\Omega}_i$. The angle preservation property of the conformal map ensures that the parameter lines on~$\hat{\Omega}_i$ are smoothly connected across element edges, see Figure~\ref{fig:mbc_3}. For the sequence of transformations from the unit square to the smoothly parameterised  patch $\hat{\Omega}_i$ we abstractly write
\begin{equation} \label{eq:conformalTotal}
	\vec{\xi}_i = \vec S_i (\vec \eta) \, , 
\end{equation}
where the non-linear mapping~$\vec S_i$ is comprised  of the conformal transformation~$z^{4/v}$ and a rotation.  It it worth emphasising that it is straightforward to compute the inverse and derivatives of the mapping~\eqref{eq:conformalTotal}. Moreover, $\vec S_i$ depends only on the valence of the considered  patch, but not on the vertex coordinates. See Figure~\ref{fig:mbc_4} for conformal parameterisation of one-ring patches with  valences $v \in \{ 3, 4, 6 \} $.
The functions~$\hat w_i $  for computing the normalised blending functions~$w_i$ are also assembled from smooth functions defined on unit squares. In our
computations~$\hat w_i$ is one quarter of a tensor-product b-spline and covers
the entire unit square. Figure~\ref{fig:mbc_5} shows the procedure for
constructing the blending function on a valence five patch. It can be seen 
that the b-spline has its maximum at the corner which maps to the central
vertex of the one-ring patch~$\hat \Omega_i$. To normalise the functions~$\hat
w_i$,  for a given point on patch~$\hat{\Omega}_i$ the corresponding point on
an overlapping patch~$\hat{\Omega}_j$ is computed with transition functions, that is, $\vec
\xi_j = ( \vec S_j \circ \vec S_i^{-1}) (\vec \xi_i)$,
cf.~\eqref{eq:conformalTotal}. The blending function is then first mapped onto
the conformal wedge and then appropriately rotated to construct the blending
function on the patch. In Figure~\ref{fig:mbc_6} the normalised blending
functions for one-ring patches with $v \in \{3,4,6 \}$ are shown. 

In \cite{ying2004simple} the weight functions~$w_i $ are chosen such that they
 are constant  in a small neighbourhood of width~$\delta$ close to the
 unit square boundaries.  This is motivated by the need to circumvent the
 singularity of the conformal map at the central vertex of the one-ring patch.
 Our numerical experiments indicate  that the finite element solutions
 are insensitive to the choice of~$\delta$ so that~$\delta=0$ is chosen. Note
 that for evaluating the finite element integrals the surface is only evaluated
 at quadrature points, which are usually away from the vertices.
 Alternatively, the one-ring patches can also be parameterised with the
 characteristic map of Catmull-Clark subdivision surfaces, which is $C^2$ smooth
 and does not have a singularity at the extraordinary vertex,
 see~\cite{antonelli2013subdivision,levin2006modified} for details.

We now proceed to the construction of a smooth surface, i.e. a two-dimensional
smooth manifold, for a given coarse control mesh. The overall approach is very
similar to the one-dimensional case introduced in the previous section. We
consider again a point~$P(x^1, x^2, x^3)$ on the manifold with a preimage with
the coordinates~$\vec \xi_i$ on the planar patch~$\hat{\Omega}_i$. Each of the
coordinate components of~$P$ are interpolated with the partition of unity
method. We write, for instance, for the~$x^1$ coordinate
\begin{equation} \label{eq:coordinatePUM2d}
	x^1 (\vec \xi_i) = \sum_j w_j(\vec \xi_i) \left (  \vec \alpha_j^1 \cdot \vec p (\vec \xi_i)  \right ) \, ,
\end{equation}
where~$\vec p(\vec \xi)$ is a vector containing the components of a polynomial basis and~$\vec \alpha_j^1$ are the unknown coefficients. Next, the coefficients~$\vec \alpha_j^1 $ are expressed in dependence of the vertex coordinates of the control mesh. The dimension of the polynomial basis has to be equal or smaller than the number of vertices in a  patch~$\hat{\Omega}_j$. Because on unstructured meshes the valence of vertices is not fixed, we choose a polynomial basis based on the one-ring with the smallest valence in the mesh. If a higher degree polynomial is desired, two or more rings of elements can be considered as patches~$\hat{\Omega}_j$. 
 When the dimension of the polynomial basis is smaller than the number of
 vertices in the patch~$\hat{\Omega}_j$, the polynomial coefficients~$\vec
 \alpha_j^1 $ are determined with a local least-squares projection on each
 patch. Neglecting for the moment the patch index, and denoting the parametric
 coordinates of the~$2v+1$ vertices on patch~$\hat \Omega_i$ as~$\vec \xi_I$ and  the
 corresponding  nodal control mesh coordinates as~$x_I^1$, the least-squares fit on~$\hat \Omega_j$ reads
\begin{equation} \label{eq:lsqr}
	 \left ( \sum_{I=1}^{2v+1} \vec p (\vec \xi_I) \vec p^\trans (\vec \xi_I)\right ) \vec \alpha^1  = \sum_{I=1}^{2v+1} \vec p (\vec \xi_I) x_I^1 \, .
\end{equation}
Note for a specific polynomial basis the matrix on the left hand side depends only on the valence of the one-ring and can be precomputed and stored. We abbreviate the least-squares projection~\eqref{eq:lsqr} with 
\begin{equation} \label{eq:lsqrProject}
 	\vec \alpha^1 = \vec A \vec x^1 \, ,
\end{equation}
where $\vec A$ is the projection matrix and the vector $\vec x^1$ contains the $x^1$ coordinates of all the vertices in the one-ring, i.e. on the planar patch. 

Finally, by making use of~\eqref{eq:lsqrProject} we can write the partition of unity interpolation~\eqref{eq:coordinatePUM2d} in dependence of the vertex coordinates of the control mesh
\begin{equation}  \label{eq:coordinatePUM2}
\begin{aligned} 
	x^1 (\vec \xi_i) &= \sum_j w_j(\vec \xi_i) \left (  \vec p (\vec \xi_i) \cdot \vec A \vec x_j^1  \right ) 
								= \sum_j w_j(\vec \xi_i) \sum_{I=1}^{2v+1} \left (   \vec p (\vec \xi_i ) \cdot \vec A \right )_I x_{j,I}^1 \\
								&= \sum_{I=1}^{2v+1} \sum_j  \left ( w_j (\xi_i) \left (  \vec p (\vec \xi_i) \cdot \vec A \right )_I  \right ) x_{j,I}^1 =  \sum_{I=1}^{2v+1} N_I (\vec \xi_i) x_{j,I}^1  \, .								
\end{aligned}
\end{equation} 
This gives rise to the definition of basis functions~$N_I(\xi_i)$, where~$I$ is the vertex id  on patch~$\hat \Omega_j$.  In conventional finite element implementations usually system matrices and vectors are evaluated by iterating over the elements in the control mesh. Moreover, during numerical integration the basis function values at pretabulated points in the integration  element are needed. In manifold finite elements the unit square in the ~$\eta^1 -\eta^2$ coordinate system is chosen as the integration element. Hence, the basis functions~$ N_I (\vec \xi_i)$ have to be evaluated in a given element and integration point~$\vec \eta = (\eta^1, \eta^2)$. 

\begin{figure}[!htbp]
  \centering
  \includegraphics[scale=1]{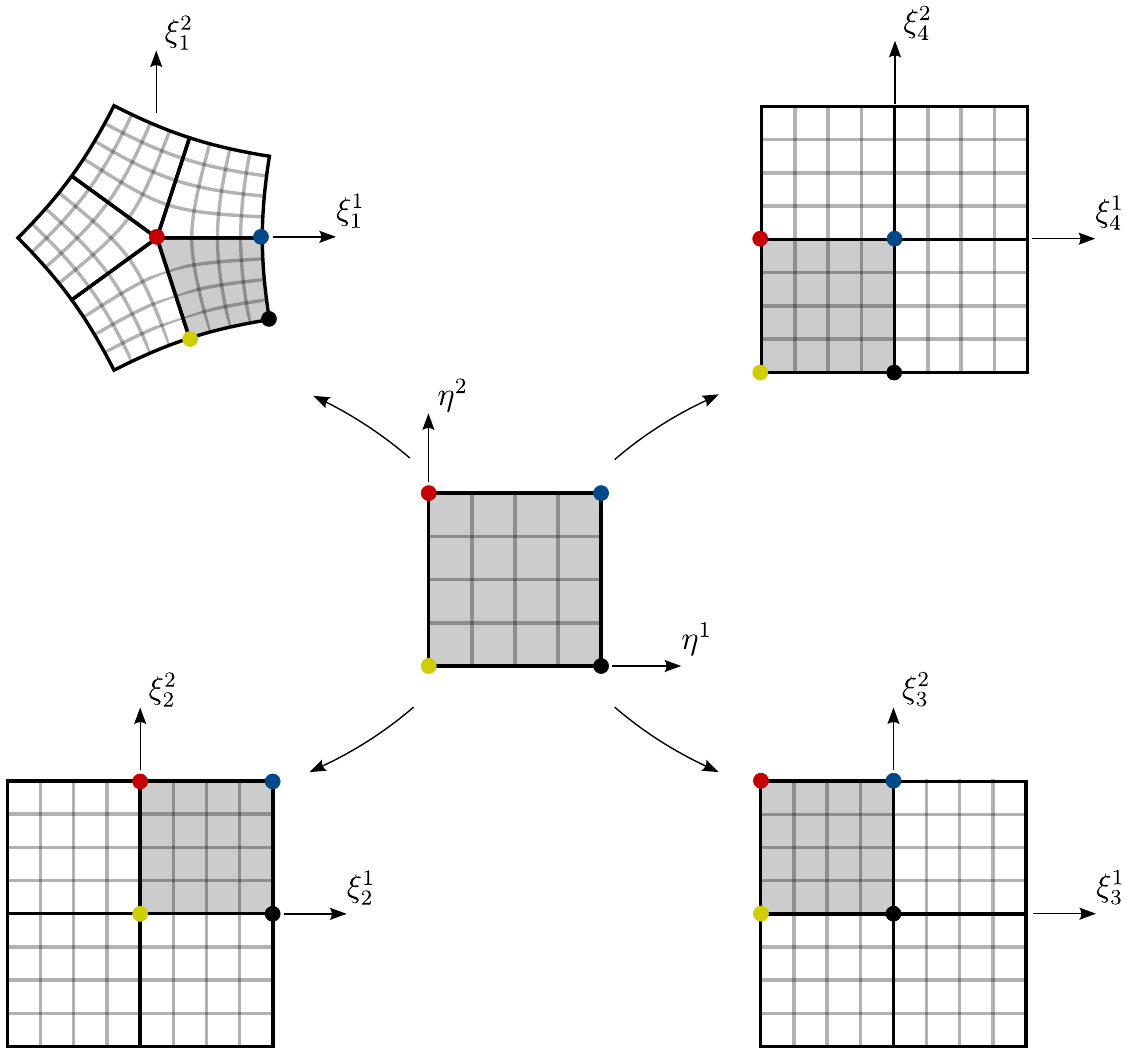}
  \caption{One integration element (centre) and its overlapping four conformally parameterised planar patches. The integration element represents the highlighted blue element in the hand control mesh shown in Figure~\ref{fig:introducEx} (left). 
  \label{fig:mbc_8}}
\end{figure}

\begin{figure}
	\centering
	\begin{minipage}{.32\textwidth}
		\centering
		\includegraphics[scale=0.2]{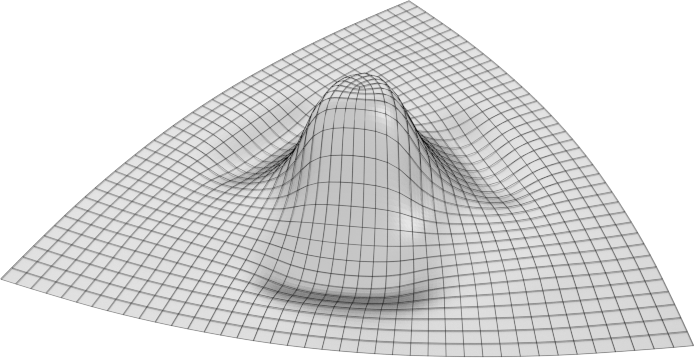}
	\end{minipage}
	\begin{minipage}{.32\textwidth}
		\centering
		\includegraphics[scale=0.2]{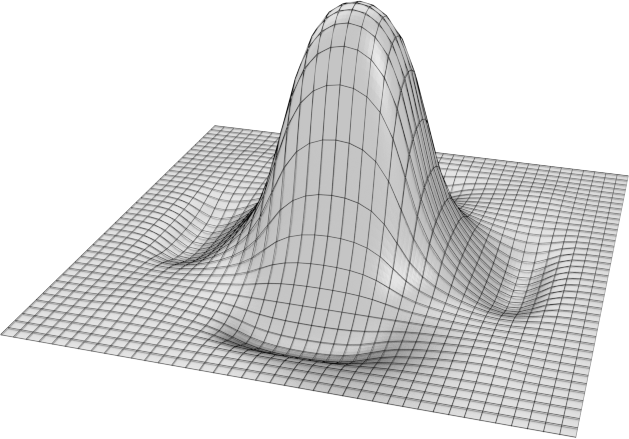}
	\end{minipage}
	\hspace*{0.005\linewidth} 	
	\begin{minipage}{.32\textwidth}
		\centering
		\includegraphics[scale=0.2]{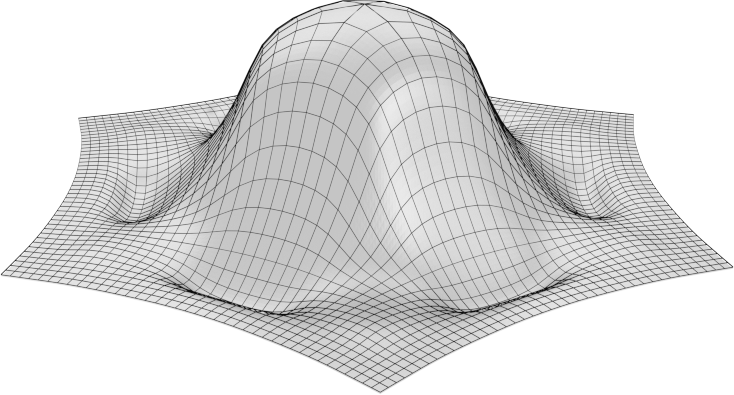}
	\end{minipage}
	\caption{Basis functions for vertices with valence three (left), four (middle)
	and six (right). The shown fine meshes  represent the parameter lines with either $\eta^1=\text{const.}$ or $\eta^2 = \text{const.}$.}
	\label{fig:mbc_9}
\end{figure}

We briefly consider the hand geometry shown in Figure~\ref{fig:introducEx} for
illustrating the process of evaluating the basis functions; see
\ref{sec:implementation} for more details. On the control mesh in
Figure~\ref{fig:introducEx} (left) one element is highlighted and the four
one-rings belonging to its four vertices are indicated. One of the element's vertices has
valence~$v=5$ and the other three have valence~$v=4$. The interpolation within
the highlighted element depends on the eighteen vertices in the union of the 
four one-rings. In Figure~\ref{fig:mbc_8} the four patches used for partition of
unity construction are shown. The unit square in the centre represents the
integration element, equivalent to the parent element in isoparametric finite
elements.  The conformal mapping~$\vec S_i (\vec \eta)$ of the unit square to
the four planar one-ring patches~$\hat \Omega_i$ depends on the valence  of the respective one-ring. For the three patches with valence four, the mapping is (up to some rotations) essentially an identity map. With the transition maps implied by the mappings~$\vec S_i$ it is straightforward to compute the basis functions defined in~\eqref{eq:coordinatePUM2}  for a given integration point~$\vec \eta_Q = (\eta^1_Q, \eta^2_Q)$. In this specific example the minimum number of vertices in a patch is nine so that  the polynomial basis $\vec p (\vec \xi_i)$  can be chosen either as a bilinear or  biquadratic Lagrangian basis.  Note that on patches with valence $v=4$ and a biquadratic Lagrangian basis the least-squares projection matrix $\vec A$ is an identity matrix. Moreover, as also can be deduced from Figure~\ref{fig:mbc_8} each basis function~$N_I$ has a support size consisting of two rings of elements around its associated vertex. In Figure~\ref{fig:mbc_9} basis functions for vertices with valence $v \in \{ 3,4,6 \}$ are shown.

\section{Examples}\label{sec:examples}  
%
We consider second and fourth-order partial differential equations to demonstrate the accuracy and convergence of the introduced manifold basis functions when used in finite element analysis.  For manifold construction, around each vertex patches consisting of either one or two-rings of elements are considered.  As blending functions we use either normalised linear, quadratic or cubic b-splines. The  element integrals are evaluated with\mbox{ $9  \times 9 $} Gauss integration points in all examples.  This high number of integration points has been chosen in order to minimise the effect of integration errors.  In convergence  studies the control meshes are refined with the Catmull-Clark subdivision scheme~\cite{Catmull:1978aa}.  The number of extraordinary vertices in a mesh remains constant because the new vertices introduced during the refinement are all ordinary. 

Surfaces with boundaries require  modified charts for  manifold constructions,
such as those introduced in~\cite{tosun2011manifold}. The specialised  treatment
of elements close to the boundary can be avoided by introducing ghost elements just outside the domain. This is achieved by reflecting  sufficient number of internal elements and vertices along the boundary.  Furthermore, the manifold basis functions are non-interpolating at the boundaries. Therefore, we use the penalty method for applying Dirichlet boundary conditions.

%
\subsection{Two-dimensional Poisson problems \label{sec:poisson2D}}
%

%
\subsubsection{Square domain with a structured mesh \label{sec:structured2D}}
\begin{figure}
\centering
	\begin{minipage}{.48\textwidth}
		\centering
		\includegraphics[scale=1.4]{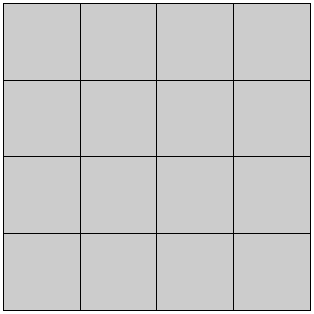}
	\end{minipage}
 	\hfill
	\begin{minipage}{.48\textwidth}
		\centering
		\includegraphics[scale=1.4]{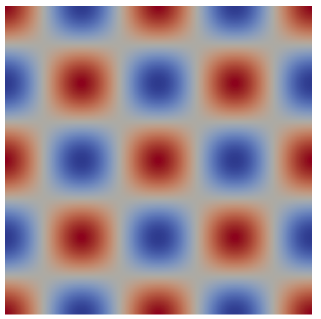}
	\end{minipage}
  \caption{Structured coarse control mesh (left) and contour plot of the solution \mbox{$u^h$} (right, using a finer control mesh).}
  \label{fig:ex1_1}
\end{figure}

\begin{figure}
  \centering
  \includegraphics[scale=0.5]{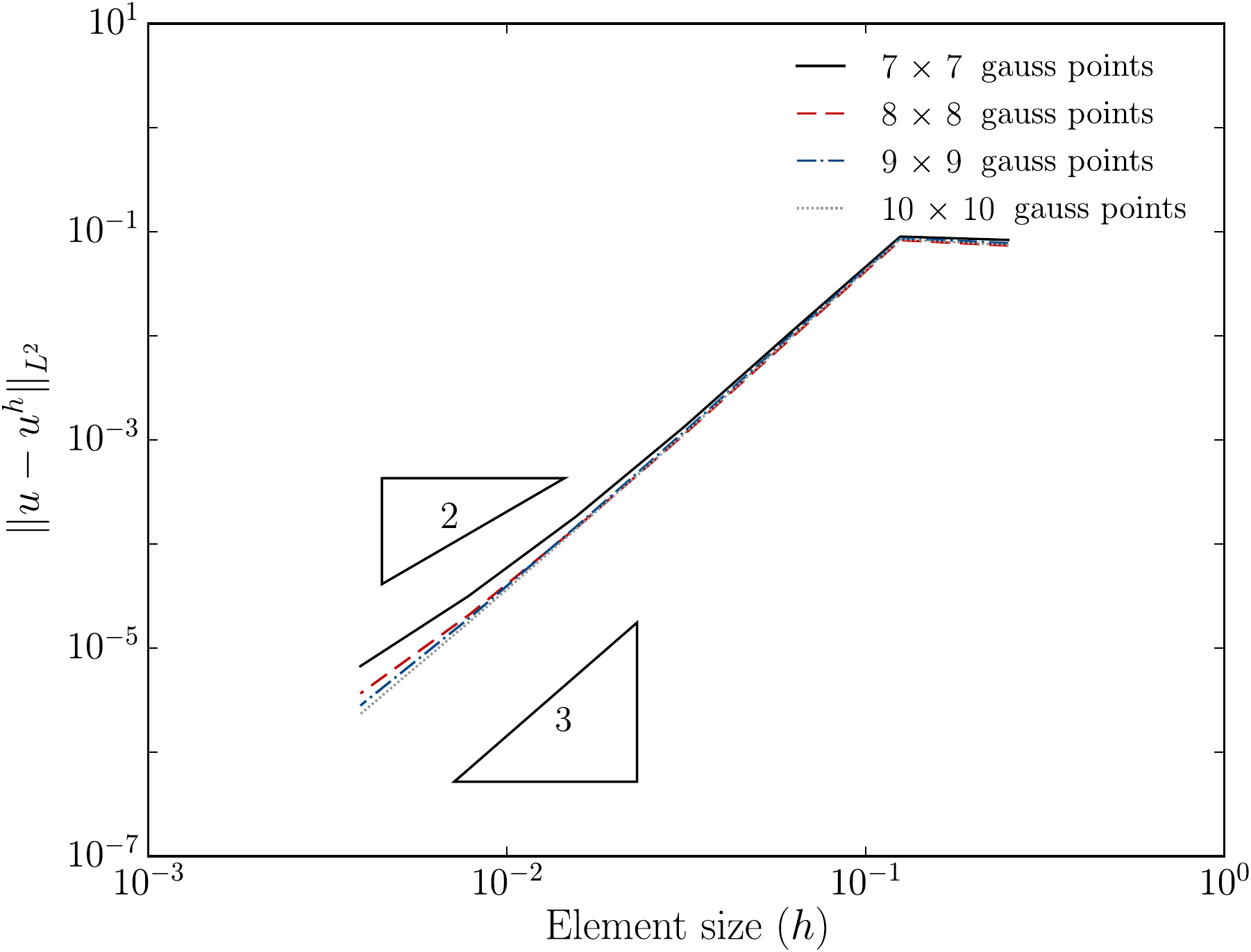}
  \caption{Square domain with a structured mesh. Influence of number of Gauss integration points on $L^2$ norm error. The local polynomial basis~$\vec p$ is quadratic and the blending functions~$w_i$ are normalised cubic b-splines. }
  \label{fig:ex1_3} 
\end{figure}

\begin{figure}
\centering
	\begin{minipage}{.49\textwidth}
		\centering
		\includegraphics[scale=0.44]{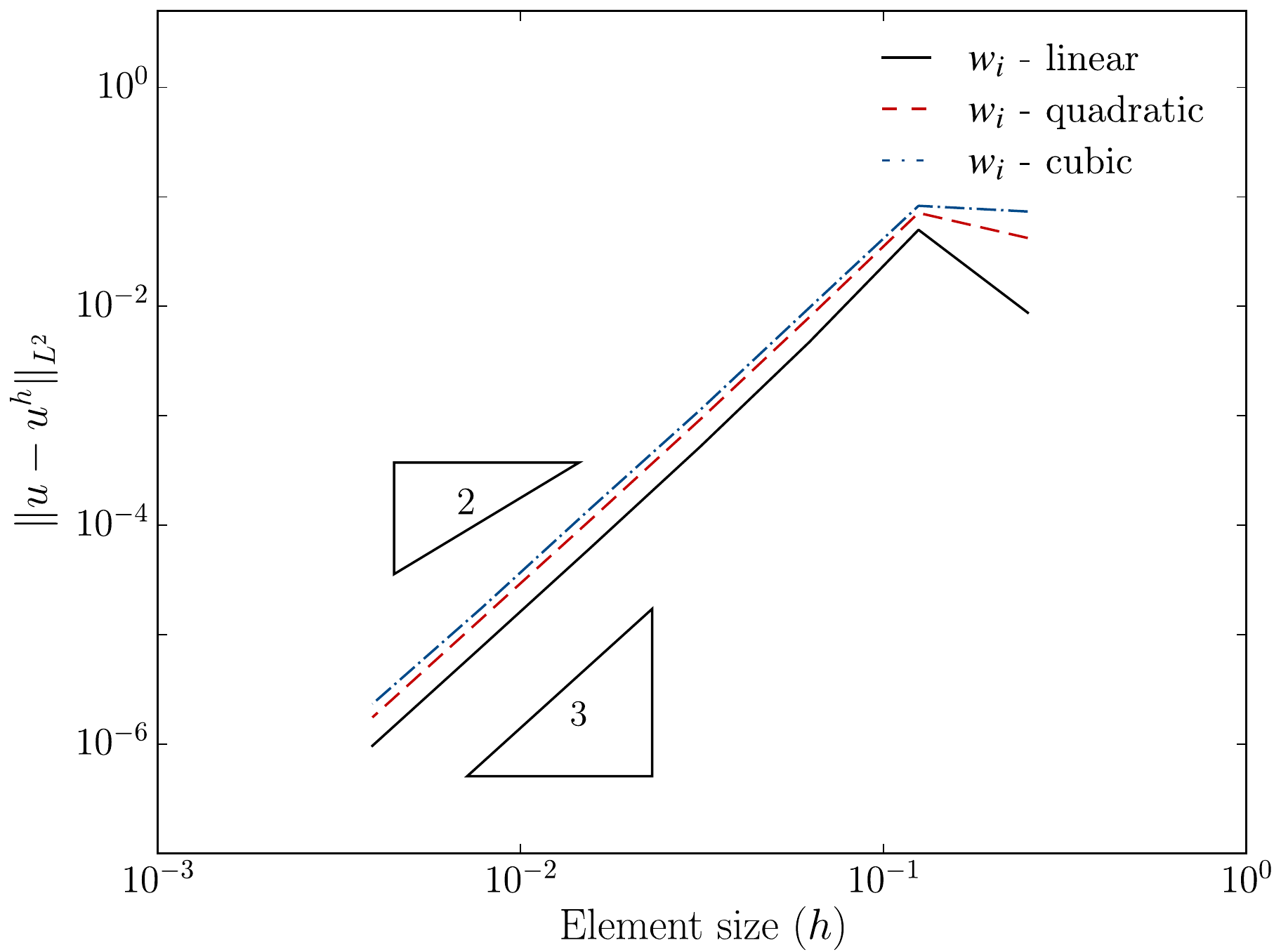}
	\end{minipage}
 	\hspace*{\fill}
	\begin{minipage}{.49\textwidth}
		\centering
		\includegraphics[scale=0.44]{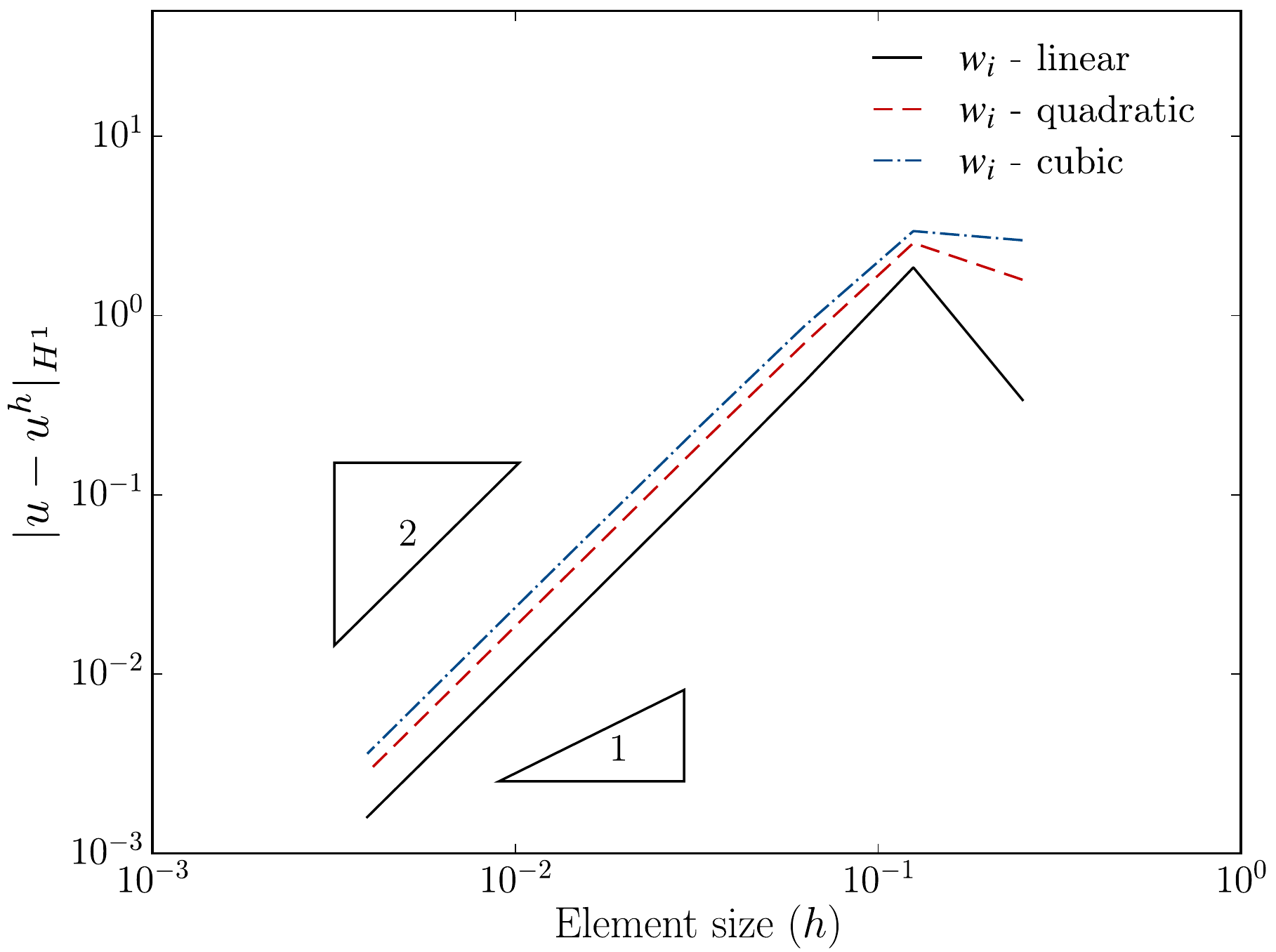}
	\end{minipage}
  \caption{Square domain with a structured mesh. Convergence of the $L^2$ norm (left) and \mbox{$H^1$} (right)
	semi-norm errors. The local polynomial basis~$\vec p$ is quadratic and the blending functions~$w_i$ are normalised linear, quadratic or cubic b-splines. 
} \label{fig:ex1_4}
\end{figure}

%
As an introductory example, we solve the Poisson-Dirichlet problem \mbox{$ -\nabla \cdot \nabla \vec{u} = q$} on the domain \mbox{$\Omega = (0, 1) \times (0, 1)$}, discretised with a Cartesian grid, see Figure~\ref{fig:ex1_1}.  The loading $q$ is chosen such that the analytical solution is
\begin{equation}\label{eq:coscos}
   u(\vec{\mathrm x}) = \cos(4 \pi x^1) \cos(4 \pi x^2)   \, ,
\end{equation}
where the variables \mbox{$\vec x = (x^1, x^2) $} represent coordinates.
In the partition of unity construction we use a biquadratic Lagrangian basis as
the local polynomial basis~\mbox{$\vec p(\vec \xi_i)$} and as blending functions~\mbox{$w(\vec \xi_i)$} we consider normalised linear, quadratic and cubic
b-splines.  To begin with,  the number of Gauss points required for adequate
integration is determined, see Figure~\ref{fig:ex1_3}.   In
Figure~\ref{fig:ex1_3}  only normalised cubic b-splines are used. Since the
basis functions are rational, a large number of Gauss points is unavoidable,
especially for finer meshes. Furthermore, our findings indicate that blending
functions with higher continuity, in general, require more Gauss points. 
According to Figure~\ref{fig:ex1_3}, the chosen~\mbox{$9 \times 9 $} Gauss
points for all examples  in the paper seems to  provide a good trade-off between
accuracy and efficiency.
Figure \ref{fig:ex1_4} shows the~\mbox{$L^2$} norm and~\mbox{$H^1$} semi-norm of
error~\mbox{$u-u^h$} as the mesh is uniformly refined with the Catmull-Clark
scheme. For this structured mesh Catmull-Clark refinement is equivalent to
refinement by bisection. The error norms for three different blending functions,
namely  normalised  linear, quadratic and cubic splines, are shown. It can be
inferred from these convergence plots that the convergence rates are optimal
and are unaffected by the blending functions. Interestingly, the constants in the
convergence plots increase when the smoothness of the blending function is
increased.
\subsubsection{Square domain with an unstructured mesh \label{sec:unStructured2D}}
\begin{figure}
\centering
	\begin{minipage}{.49\textwidth}
		\centering
		\includegraphics[scale=1.5]{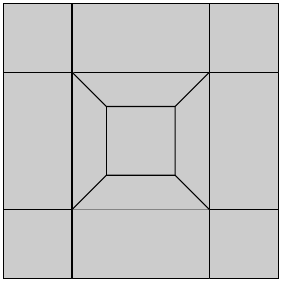}
	\end{minipage}
	\hfill
	\begin{minipage}{.49\textwidth}
		\centering
		\includegraphics[scale=1.5]{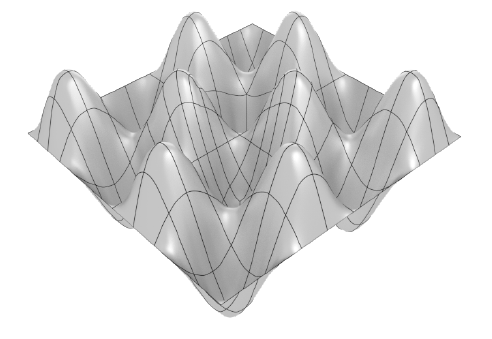}
	\end{minipage}
  \caption{ Unstructured coarse control mesh (left) and the finite element solution (right). }
  \label{fig:ex2_1}
\end{figure}

%
This example underlines the performance of manifold basis functions on unstructured meshes, and studies how the convergence rates are influenced in the presence of extraordinary vertices. To this end, the Poisson-Dirichlet problem with the analytical solution
\begin{equation}\label{eq:sinsin}
   u(\vec{\mathrm x}) = \sin(4 \pi x^1) \sin(4 \pi x^2)   
\end{equation}
is considered.  Figure \ref{fig:ex2_1} shows the unstructured mesh used in the
computations and  a representative finite element solution. The mesh lines on the displaced solution in Figure \ref{fig:ex2_1} represent the edges of the elements on the exact surface. The mesh has eight extraordinary vertices, with four vertices of valence~\mbox{$v=3$} and the other four of valence~\mbox{$v=5$}.
In the convergence studies we use as local polynomials~\mbox{$\vec p (\vec \xi_i)$} bilinear and biquadratic  Lagrangian polynomials.  As patches one- and two-rings of elements are considered. However, in one-ring patches with valence~\mbox{$v=3$} there are only seven vertices so that instead of a biquadratic Lagrangian polynomial  locally a complete quadratic polynomial has to be used.  In all cases the blending functions~\mbox{$w(\vec \xi_i)$} are normalised cubic b-splines.

\begin{figure}
\centering
	\begin{minipage}{.49\textwidth}
		\centering
		\includegraphics[scale=0.44]{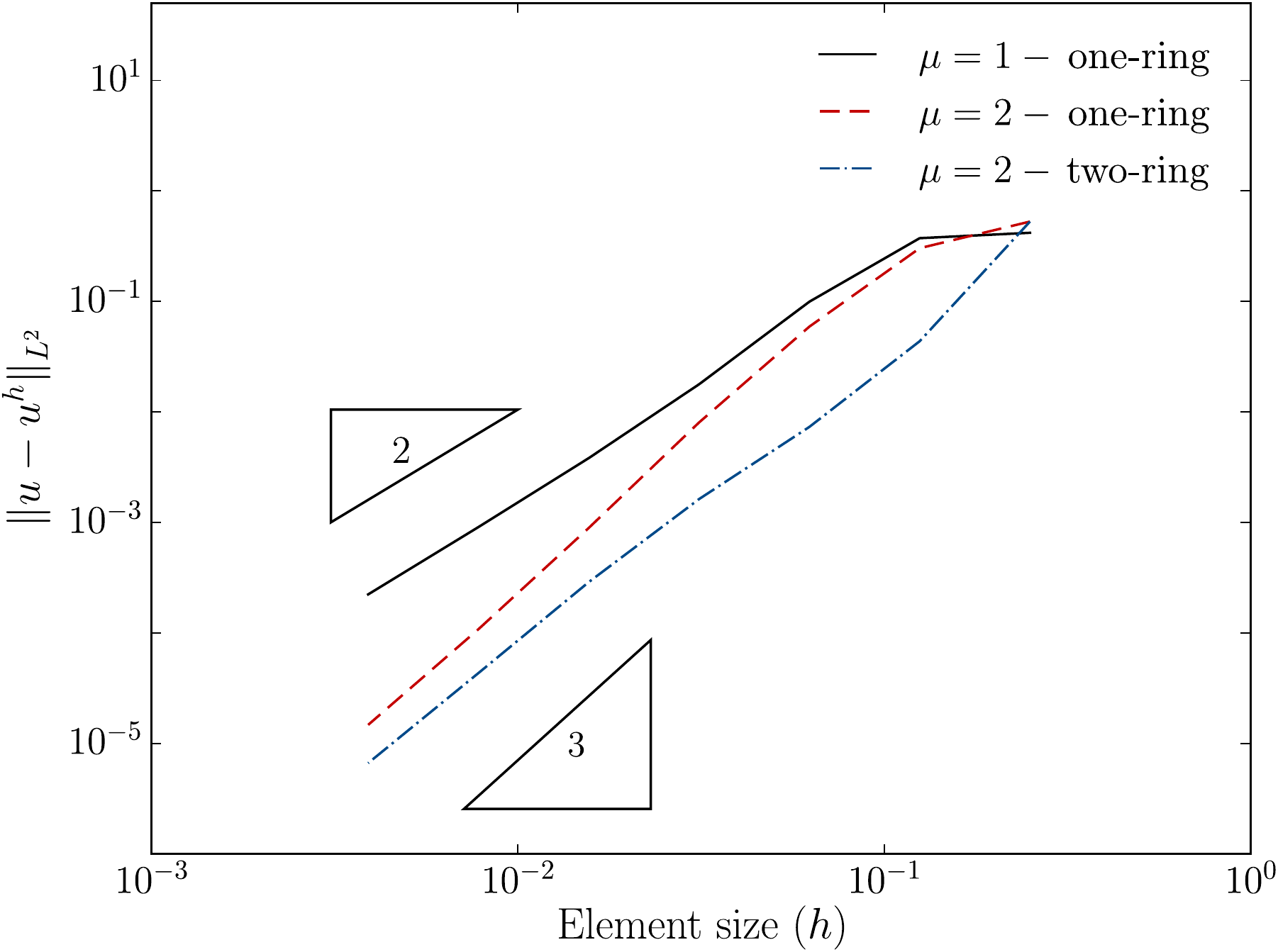}
	\end{minipage}
	\hfill
	\begin{minipage}{.49\textwidth}
		\centering
		\includegraphics[scale=0.44]{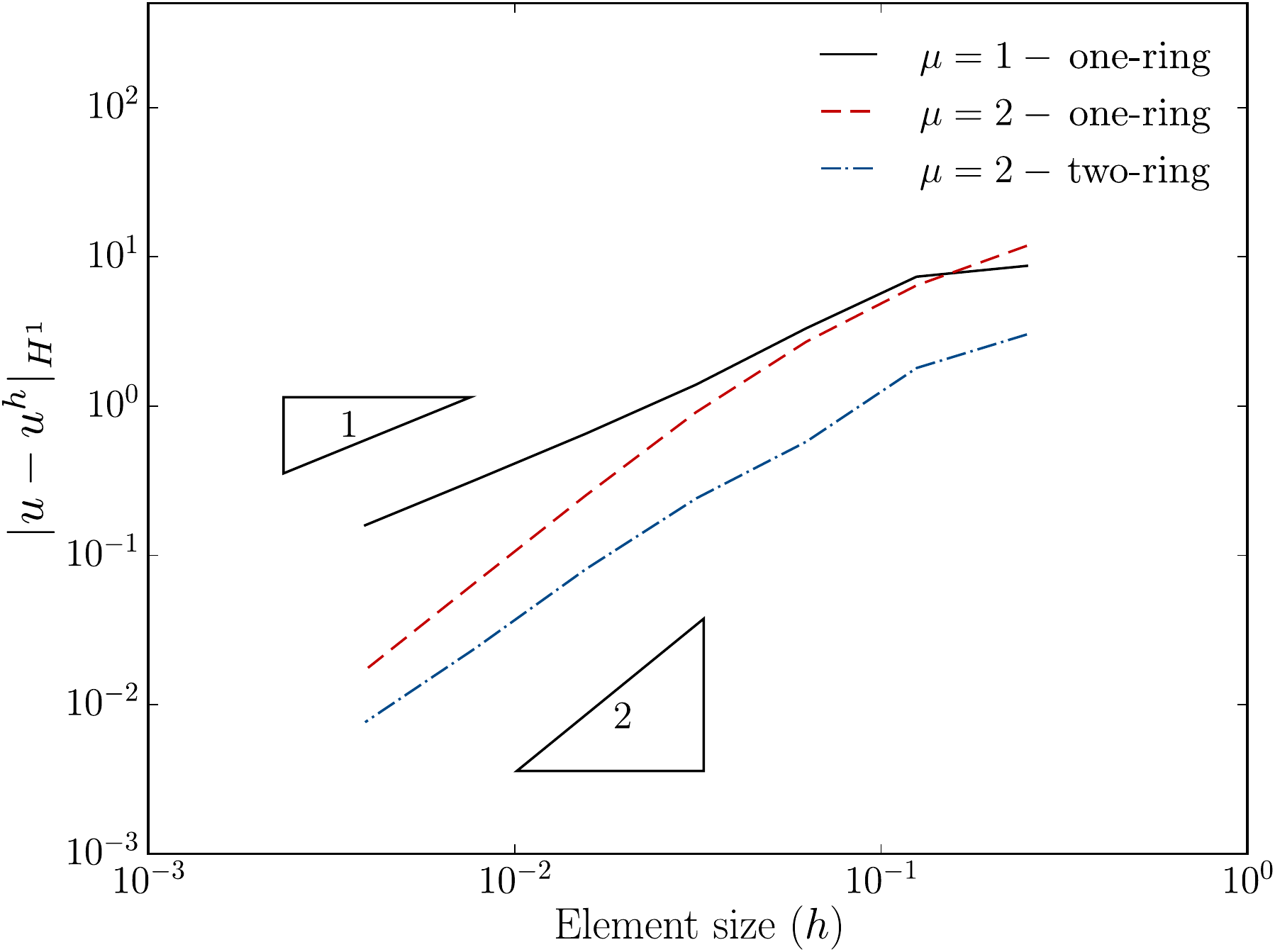}
	\end{minipage}
  	\caption{Square domain with an unstructured mesh. Convergence of the $L^2$norm (left) and $H^1$ (right) 
	semi-norm errors.  The local polynomial basis~$\vec p$ is either  linear ($\mu=1$) or quadratic ($\mu=2$) and the blending functions $w_i$  are normalised  cubic b-splines.} \label{fig:ex2_3}
\end{figure}

\begin{figure}
\centering
	\begin{minipage}{.49\textwidth}
		\centering
		\includegraphics[scale=0.53]{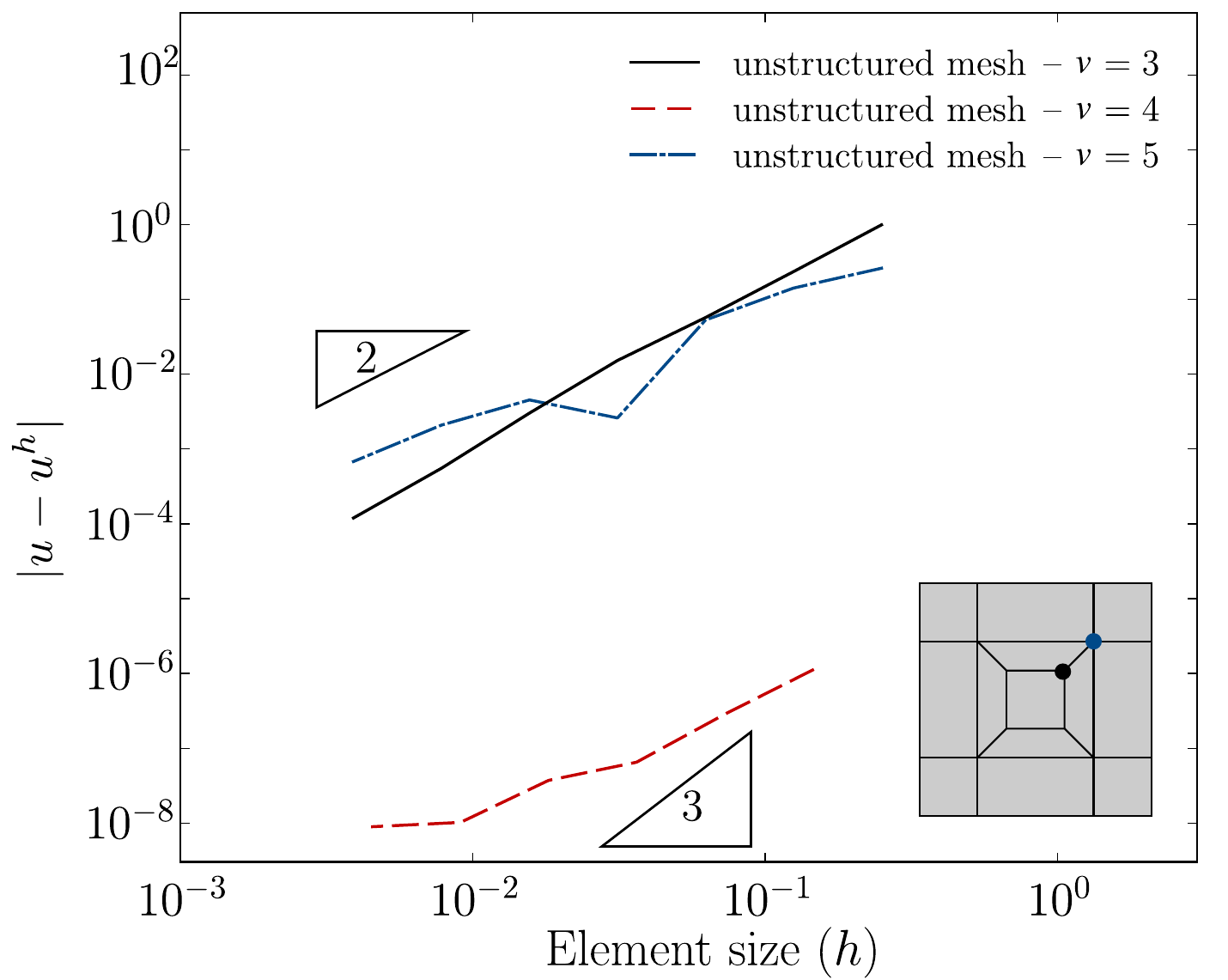}
	\end{minipage}
	\hfill
	\begin{minipage}{.49\textwidth}
		\centering
		\includegraphics[scale=0.53]{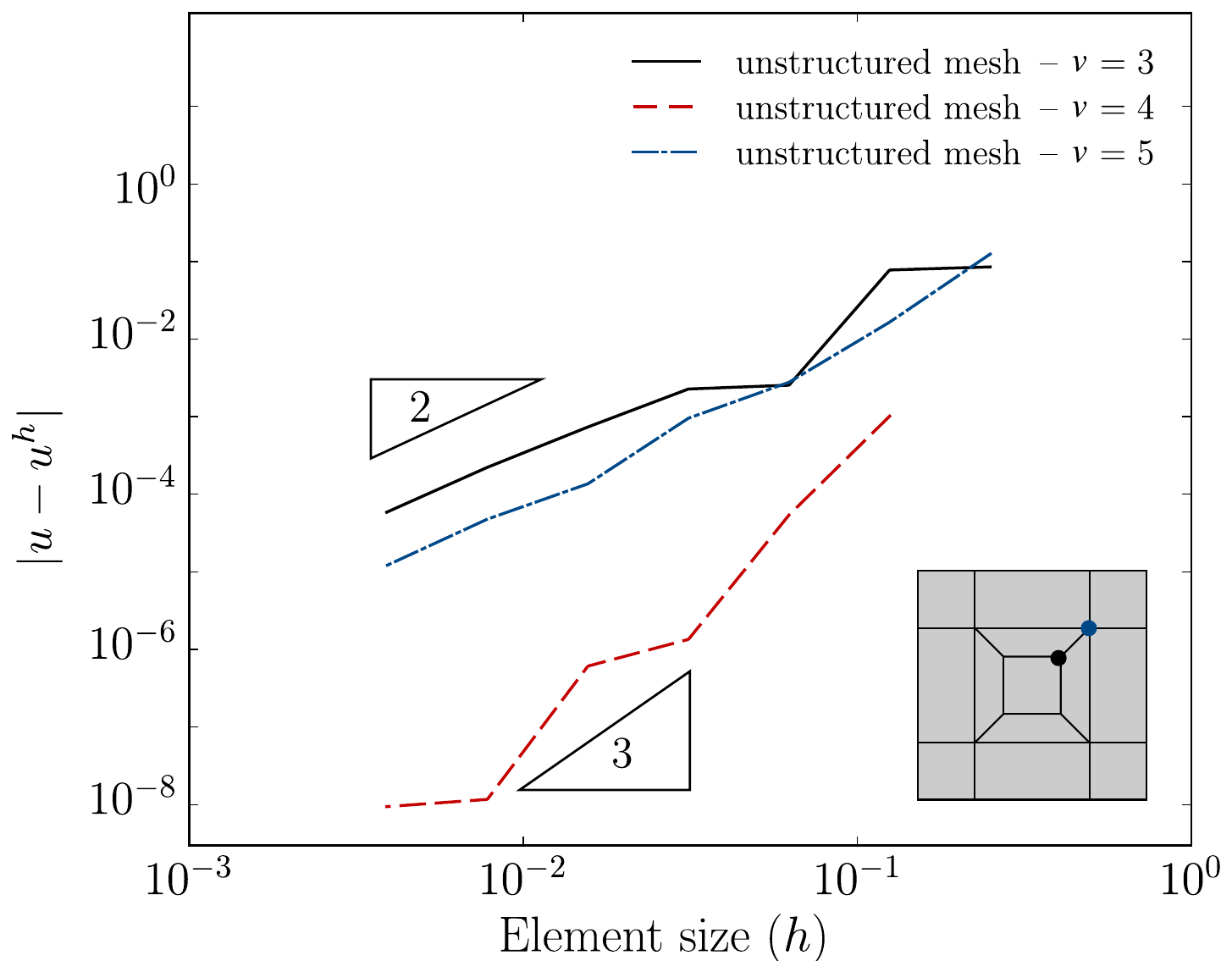}
	\end{minipage}
  	\caption{Square domain with an unstructured mesh. Pointwise convergence of the solution at three selected vertices with valences  $v \in \{3, 4, 5 \}$ for one-ring (left) and two-ring patches (right).  The local polynomial basis~$\vec p$ is quadratic and the blending functions $w_i$  are normalised  cubic b-splines.} 
	\label{fig:ex2_5}
\end{figure}

Figure \ref{fig:ex2_3} shows the~\mbox{$L^2$} norm and~\mbox{$H^1$} semi-norm 
of the error \mbox{$u-u^h$} as the mesh is successively refined
with the Catmull-Clark scheme. It can be seen that for all cases, the
convergence rates are close to optimal.  In the~\mbox{$L^2$} norm and for
quadratic polynomials ~(\mbox{$\mu = 2$})  the convergence rates for one- and
two-ring patches are approximately~$2.9$ and~$2.7$ respectively. In the
\mbox{$H^1$} semi-norm, the corresponding convergence rates are  $1.9$ and $1.7$
respectively.  We believe that the reduction of convergence rates with
increasing patch size is primarily due to the suboptimal integration of rational
polynomials.
Next, we study the pointwise convergence of the solution at selected vertices.
For this part of the studies we use one- and two-ring patches and  consider only
biquadratic Lagrangian polynomials.  As mentioned before for one-ring patches
with valence~\mbox{$v=3$} instead of the biquadratic Lagrangian a complete
quadratic has to be used. The blending functions are normalised cubic b-splines.
 Figure~\ref{fig:ex2_5} (left) shows the convergence of the error at three
 selected vertices  with valences~\mbox{$v \in \{3, 4,5 \}$} when one-ring patches are
used. Figure~\ref{fig:ex2_5} (right)  shows the corresponding plots when
two-ring patches are used. For both types of patches, the convergence rate at
extraordinary vertices is~$2$.  It is interesting to note that two-ring patches,
in general, have smaller errors than the one-ring patches.
%
%
\subsubsection{Circular domain \label{sec:circular2D}}
\begin{figure}
\centering
	\begin{minipage}{.49\textwidth}
		\centering
		\includegraphics[scale=1]{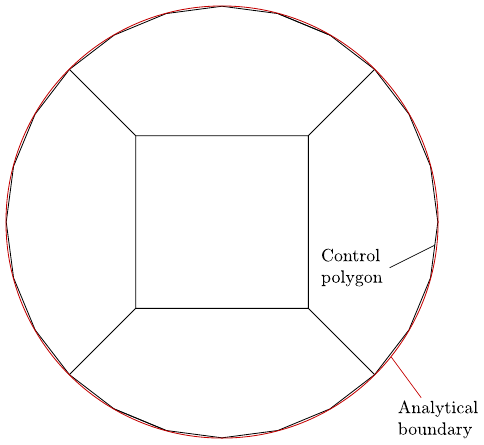}
	\end{minipage}
	\hfill
	\begin{minipage}{.49\textwidth}
		\centering
		\includegraphics[scale=0.25]{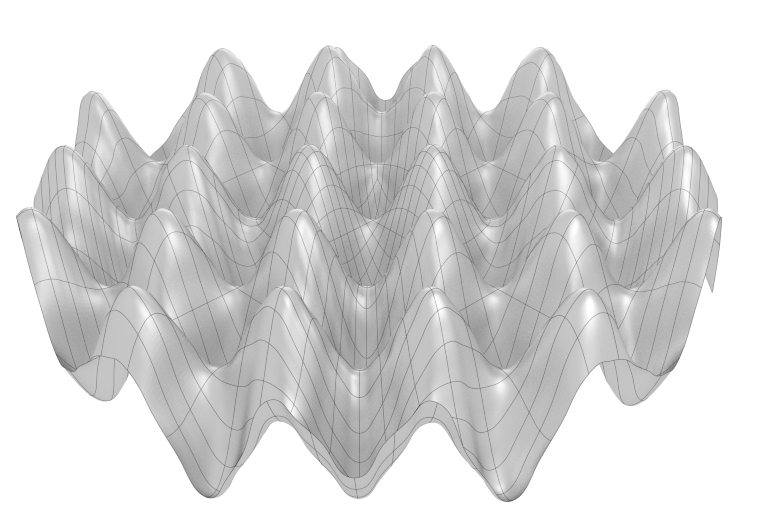}
	\end{minipage}
	\caption{Exact and approximated circular domains and the coarse
	control mesh (left) and finite element solution (right).}
	\label{fig:ex3_1}
\end{figure}

%
In this example the Poisson-Dirichlet problem with the same solution as in~\eqref{eq:sinsin} 
is solved on a circular domain with radius~\mbox{$||x||<0.5$}, see Figure~\ref{fig:ex3_1}.  The aim is
to illustrate the treatment of curved boundaries when manifold basis functions
are used. As previously mentioned, in our present implementation, ghost nodes
are introduced at the boundaries, which circumvents the use of modified  charts
close to the boundaries.  With specialised boundary charts the treatment of
boundaries would be different.

\begin{figure}[]
  \centering
  \includegraphics[scale=0.5]{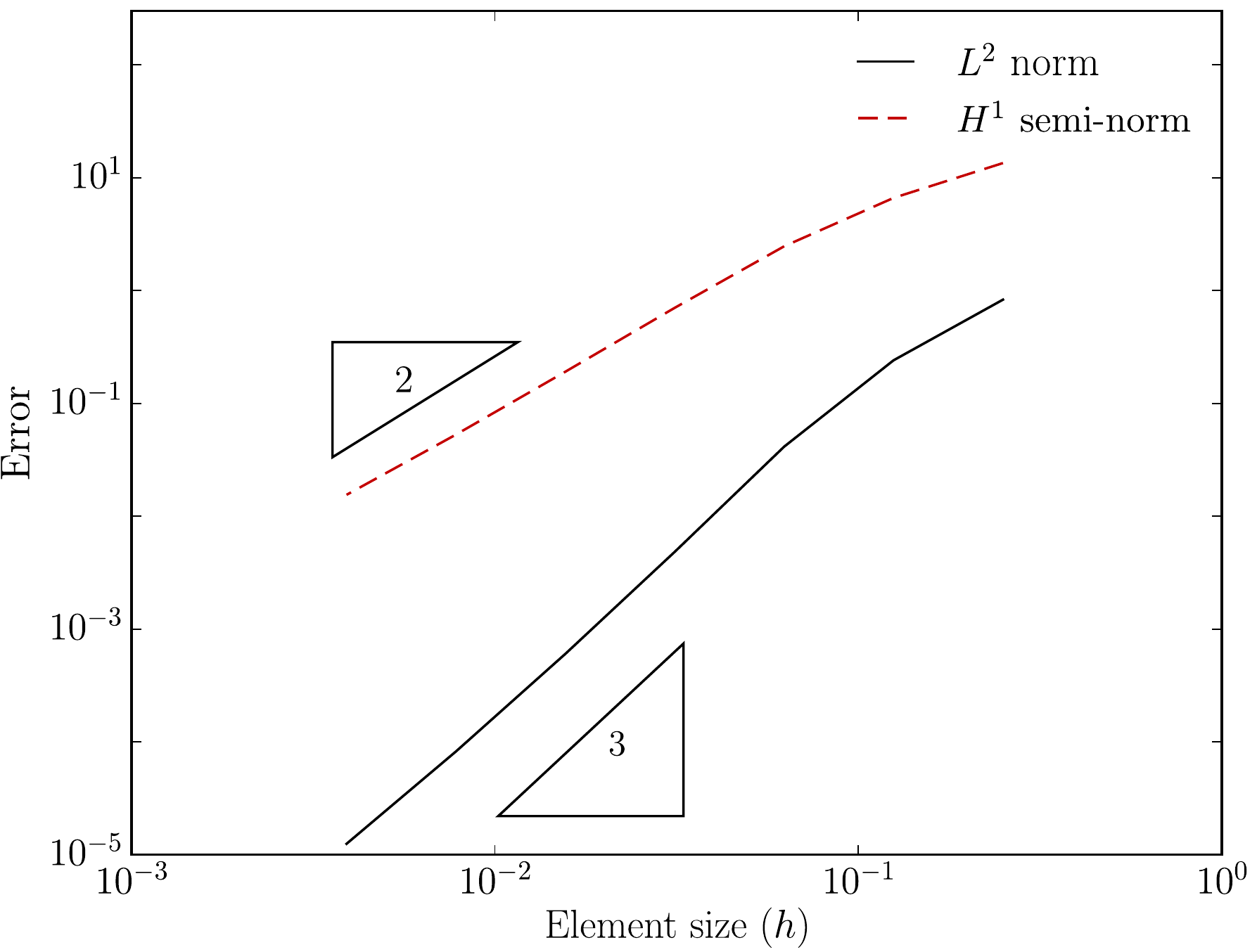}
  \caption{Circular domain. Convergence of the~$L^2$ norm and~$H^1$ semi-norm errors. 
  The local polynomial basis~$\vec p$ is quadratic and the blending functions~$w_i$  are normalised  cubic b-splines.}
  \label{fig:ex3_3}
\end{figure}

Figure~\ref{fig:ex3_1} (left) shows the coarse mesh containing four
extraordinary vertices with valence~\mbox{$v=3$}. The domain boundary is approximated by
least-squares fitting an approximate circle described by the manifold basis
functions  to the exact circle. In the least-squares problem the unknowns are
the positions of the vertices close to the boundary. This is performed as a
preprocessing step each time after the mesh is refined with Catmull-Clark
subdivision.
In the convergence study shown in Figure~\ref{fig:ex3_3} only one-ring patches
are used. The local polynomials~\mbox{$\vec p(\vec \xi_i)$} are either
biquadratic Lagrangian polynomials for patches with valence~\mbox{$v \neq 3$} or
complete quadratic polynomials for patches with \mbox{$v = 3$}.  The blending
functions~\mbox{$w(\vec \xi_i)$} are normalised cubic b-splines. As can
be inferred from Figure~\ref{fig:ex3_3}, optimal convergence rates for the~$L^2$
norm and the~$H^1$ semi-norm error are achieved.

%
\subsection{Thin-plate and thin-shell problems \label{sec:platesAndShells}}
%
The linear Kirchhoff-Love model is used for computing the thin-plate and
thin-shell problems introduced in the following. The corresponding weak form
depends on the metric and curvature tensors of the mid-surface in the reference
and deformed configurations. Due to the presence of the curvature tensor, the
basis functions have to be smooth, or more precisely in space~$H^2$. The
presented computations have been performed by replacing the subdivision basis
functions in our software~\cite{Cirak:2011aa, Long:2012aa} with manifold-basis
functions. Due to the algorithmic similarities between subdivision and manifold basis functions it is
straightforward to replace one with the other. 
Out-of-plane shear deformations have been neglected although it would
be possible to take them into account as proposed  in~\cite{Long:2012aa}.

%
\subsubsection{Simply supported square plate \label{sec:squarePlate}}
%

We consider the deformation of a simply supported square plate of unit length
subjected to an applied pressure loading~\mbox{$q =10000$}. The thickness of the plate is~\mbox{$h=0.01$}, the Young's modulus is \mbox{$E = 70 \times 10^9$}, and the Poisson ratio is~\mbox{$\nu=0.3$}.  Its  analytical solution according to \cite{Timoshenko:1959aa} is
\begin{equation}
   u(\vec{\mathrm x}) = \frac{16 q}{\pi^6 D} \sum_{i=1}^{\infty}
   \sum_{j=1}^{\infty} \frac{\sin(i \pi x^1) \sin(j \pi x^2)}{ij
   (i^2 + j^2)^2} \, ,
\end{equation}
where $D = Eh^3/(12 (1 - \nu^2)$ is the flexural rigidity.  

\begin{figure}
\centering
	\begin{minipage}{.45\textwidth}
		\centering
		\includegraphics[scale=1]{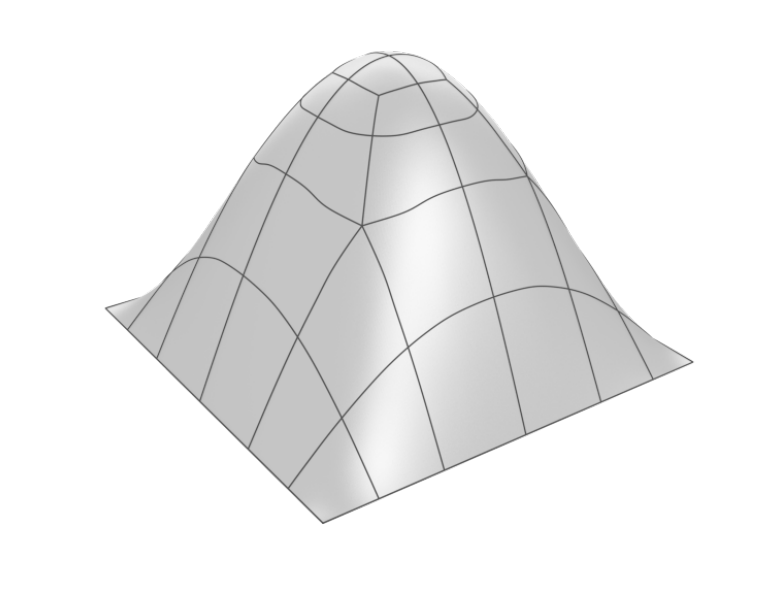}
	\end{minipage}
\hfill
	\begin{minipage}{.53\textwidth}
		\centering
		\includegraphics[scale=0.45]{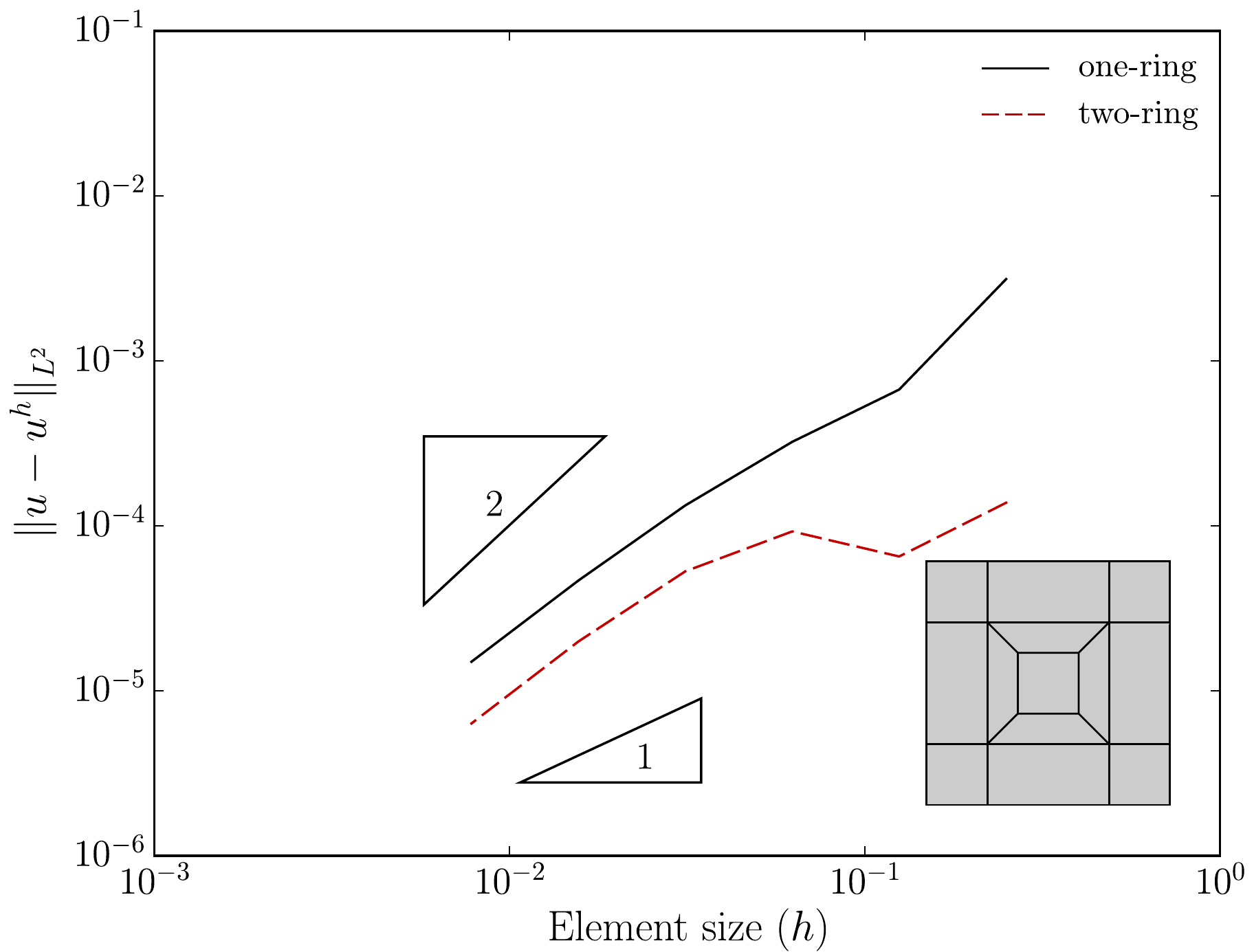}
	\end{minipage}
	\caption{Simply supported square plate. Deflected plate  (left) and convergence of the $L^2$ norm  error (right). The local polynomial basis~$\vec p$ is quadratic and the blending functions $w_i$  are normalised  cubic b-splines.}
	\label{fig:plate_1}
\end{figure}

The unstructured mesh shown in Figure \ref{fig:ex2_1}, previously used for the
Poisson-Dirichlet problem, is reused for this plate problem. Recall that the mesh
has eight extraordinary vertices, namely  four with~\mbox{$v=3$} and the other
four with~\mbox{$v=5$}.  In the convergence study, one-rings and two-rings of
elements are considered.  Except on valence~\mbox{$v=3$} one-ring patches,  the
local polynomials~\mbox{$\vec p ( \vec \xi_i)$} are biquadratic Lagrangian
polynomials. As previously mentioned,  on  one-ring patches with~$v=3$ only a
complete quadratic polynomial can be used, because there are only seven vertices
in a patch. In all cases the blending function~$w(\vec \xi_i)$  is a normalised
cubic b-spline. In  Figure \ref{fig:plate_1} (left), the deflected shape of the
simply-supported plate for a relatively coarse mesh  is shown.  Figure
\ref{fig:plate_1} (right) illustrates the convergence of the~$L^2$ norm  error
as the mesh is successively refined with the Catmull-Clark scheme. The constant in the
convergence plot decreases with increasing patch size.  For both patch sizes,
the convergence rate is approximately~$1.7$, which is slightly lower than the
optimum of~$2$. One possible reason for this is the inadequate integration of
the rational polynomials.
%

%
\subsubsection{Pinched cylinder \label{sec:cylinder}}
\begin{figure}
\centering
	\begin{minipage}{.49\textwidth}
		\centering
		\includegraphics[scale=0.42]{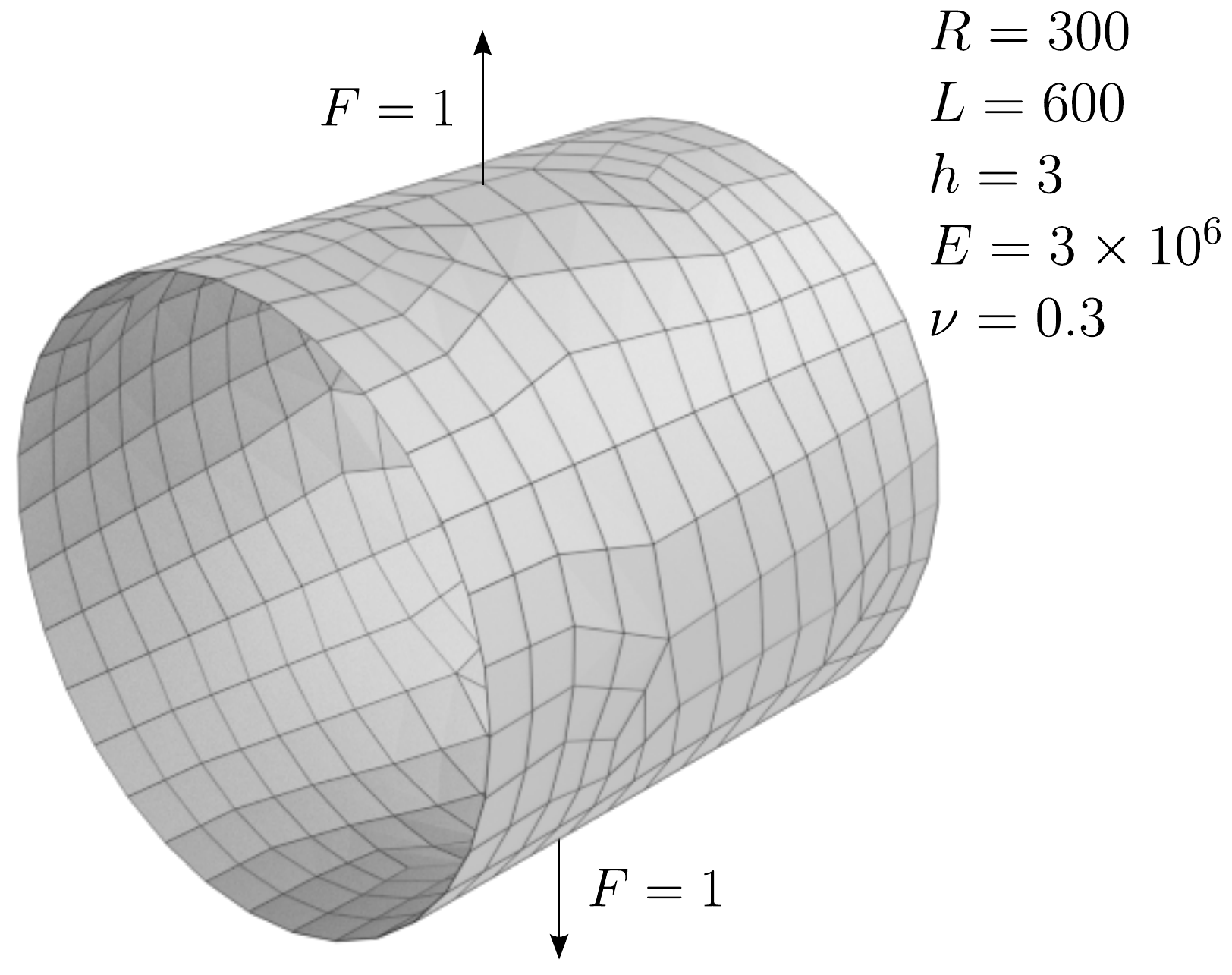}
	\end{minipage}
	\begin{minipage}{.49\textwidth}
		\centering
		\includegraphics[scale=0.43]{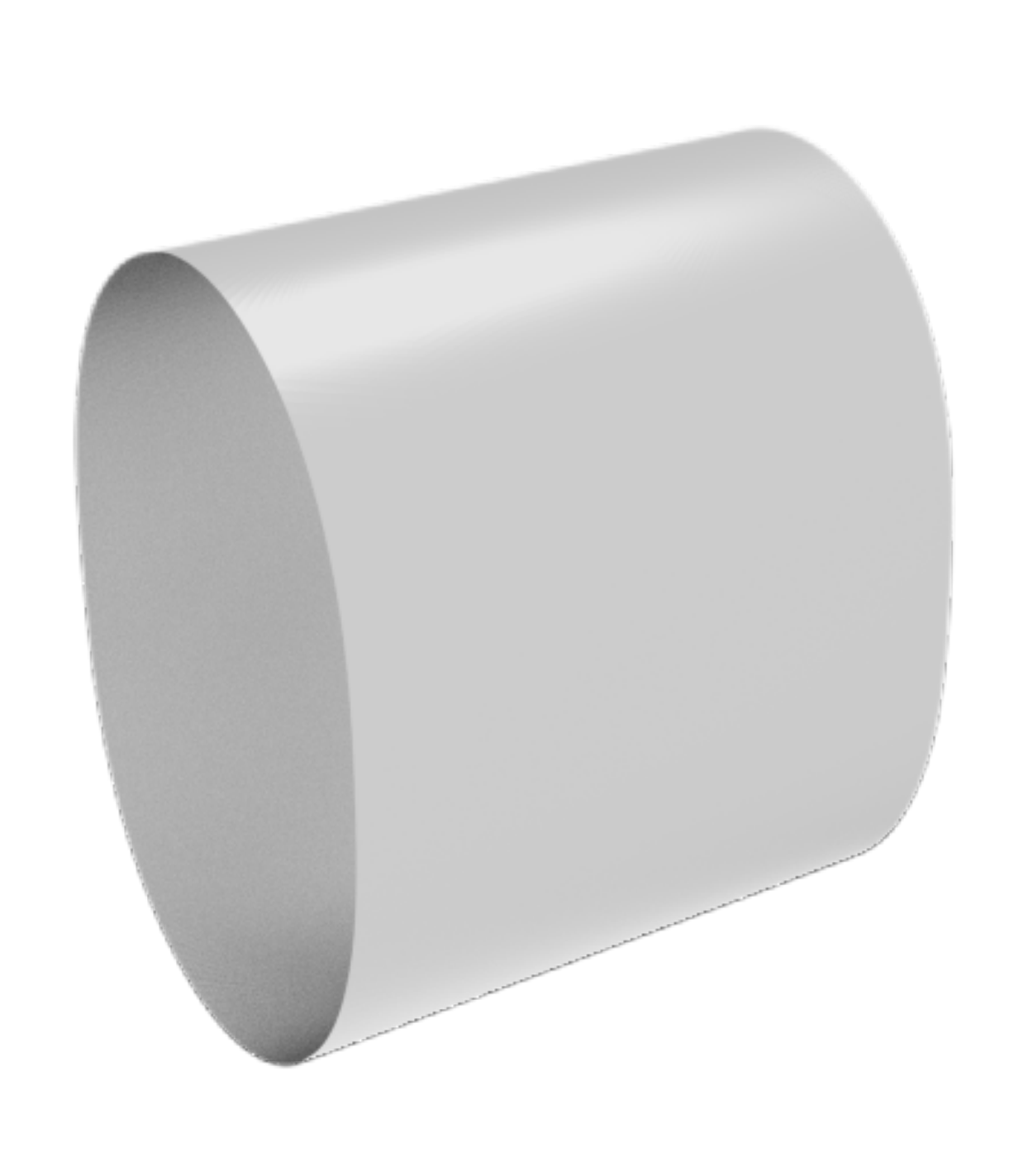}
	\end{minipage}
	\caption{Definition of the pinched cylinder problem
	(left) and deformed pinched cylinder (right).}
	\label{fig:ex5_1}
\end{figure}

\begin{figure}
\centering
	\begin{minipage}{.49\textwidth}
		\centering
		\includegraphics[scale=1.1]{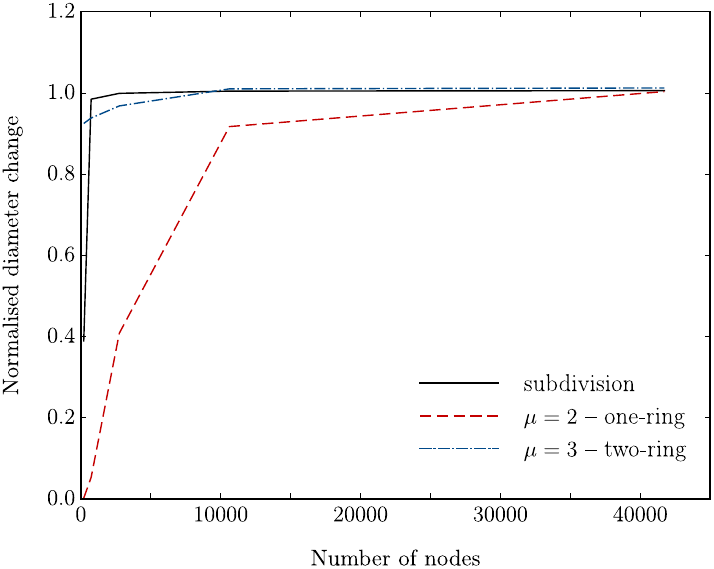}
	\end{minipage}
\hfill
	\begin{minipage}{.49\textwidth}
		\centering
		\includegraphics[scale=1.1]{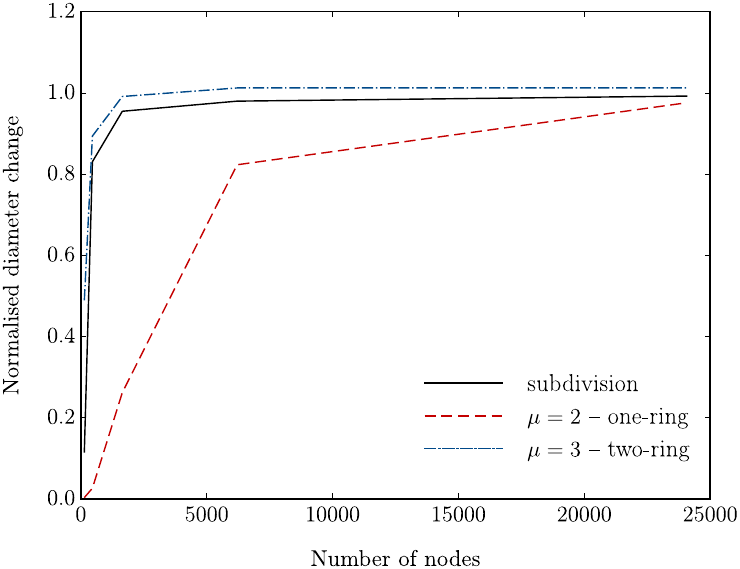}
	\end{minipage}
	\caption{Pinched cylinder. Convergence of maximum change in diameter with manifold
	and Catmull-Clark basis functions for structured (left) and unstructured
	(right) meshes.  The local polynomial basis~$\vec p$ is either quadratic ($\mu=2$) or cubic ($\mu=3$) and the blending functions $w_i$  are normalised  cubic b-splines.}
	\label{fig:ex5_3}
\end{figure}

%
The pinched cylinder is one of the benchmark examples for shell finite elements
suggested in Belytschko et al.~\cite{Belytschko:1985aa}. The unstructured coarse
control mesh, material properties and the general problem setup are shown in
Figure~\ref{fig:ex5_1} (left).  The ends of the cylinder are unconstrained and
the two diametrically opposite  forces are applied within the middle section of
the cylinder.
The deflected pinched cylinder  is shown in Figure~\ref{fig:ex5_1}
(right). In Figure~\ref{fig:ex5_3} the convergence of the normalised maximum change in diameter of the Catmull-Clark and manifold finite element solutions are compared, both structured and unstructured meshes are considered.  The change in diameter is normalised with the exact solution~$4.52 \times 10^{-4}$ for a membrane shell given in~\cite{Timoshenko:1959aa}. It can be seen that the Catmull-Clark and manifold solutions converge to very similar values. The convergence of manifold functions with cubic polynomials~($\mu = 3$) on two-ring patches is comparable with subdivision basis functions. Note that for structured meshes Catmull-Clark subdivision basis functions are identical to tensor-product cubic b-splines~\cite{Peters:2008aa}.  The slower convergence of the manifold basis functions  with quadratic polynomials~($\mu = 2$) on one-ring patches is as expected.
%
\subsubsection{Pinched hemisphere \label{sec:hemisphere}}
\begin{figure}
\centering
	\begin{minipage}{.49\textwidth}
		\centering
		\includegraphics[scale=0.45]{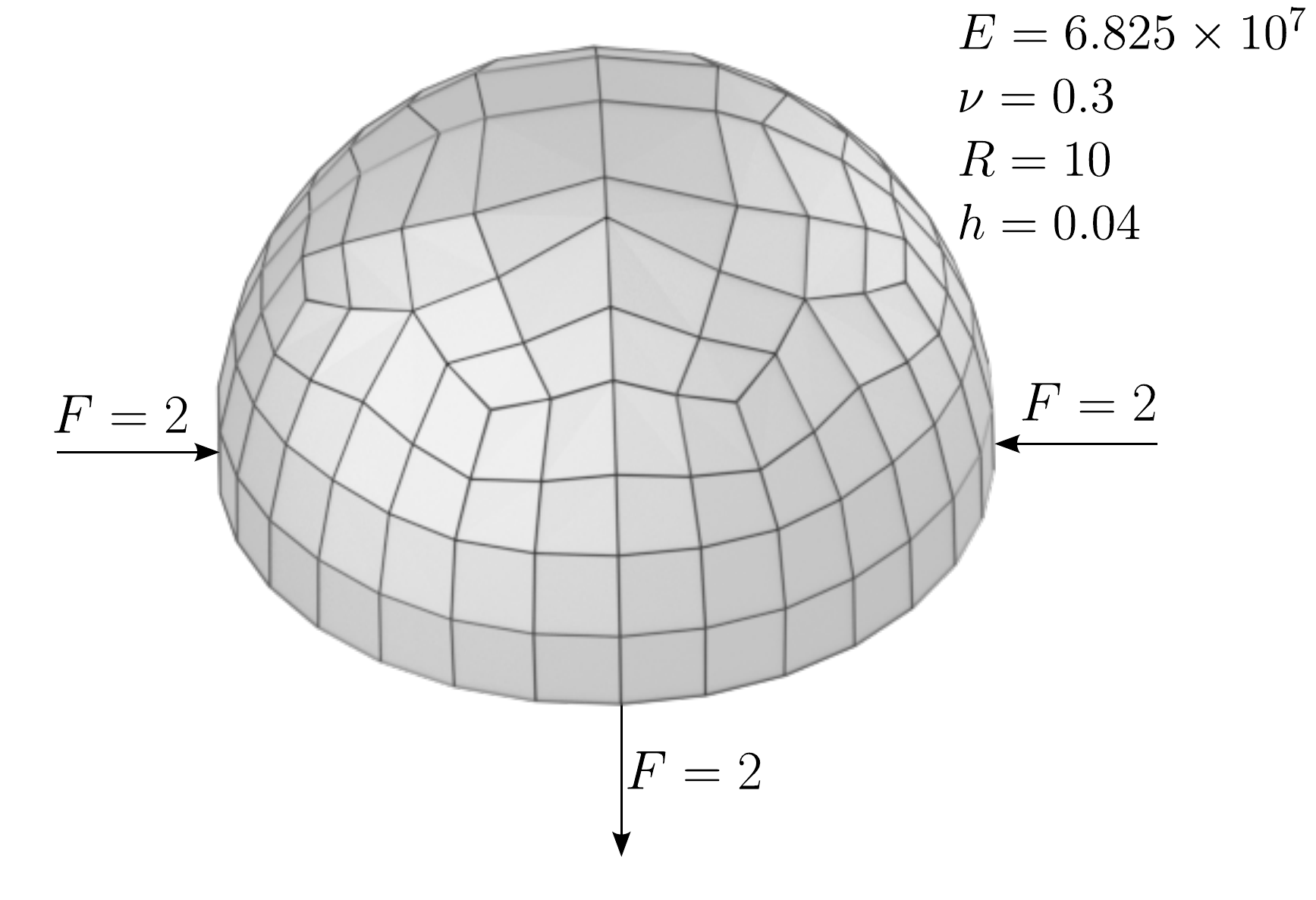}
	\end{minipage}
 	\hfill
	\begin{minipage}{.49\textwidth}
		\centering
		\includegraphics[scale=0.3]{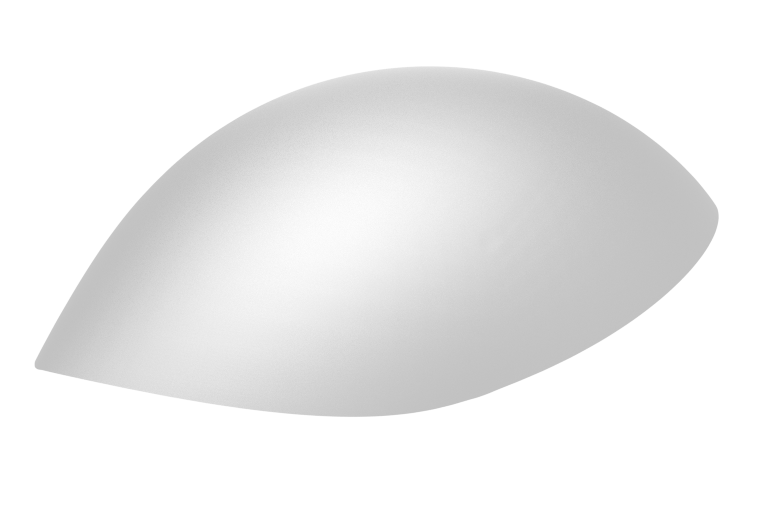}
	\end{minipage}
	\caption{Definition of the pinched  hemisphere problem
	(left) and deformed pinched hemisphere (right).}
	\label{fig:hemi_1}
\end{figure}

%
Our last example is the pinched hemisphere, which has also been suggested in
Belytschko et al.~\cite{Belytschko:1985aa} as a benchmark for shell finite
elements. The coarse control mesh, material properties and the general problem setup are shown in
Figure \ref{fig:hemi_1} (left). The edge of the hemisphere is unconstrained and
the four radial forces have alternating signs. The sum of the applied forces is zero. In the
control mesh the valence of the vertices  range between 3 and 5.

\begin{figure}[]
  \centering \includegraphics[scale=0.5]{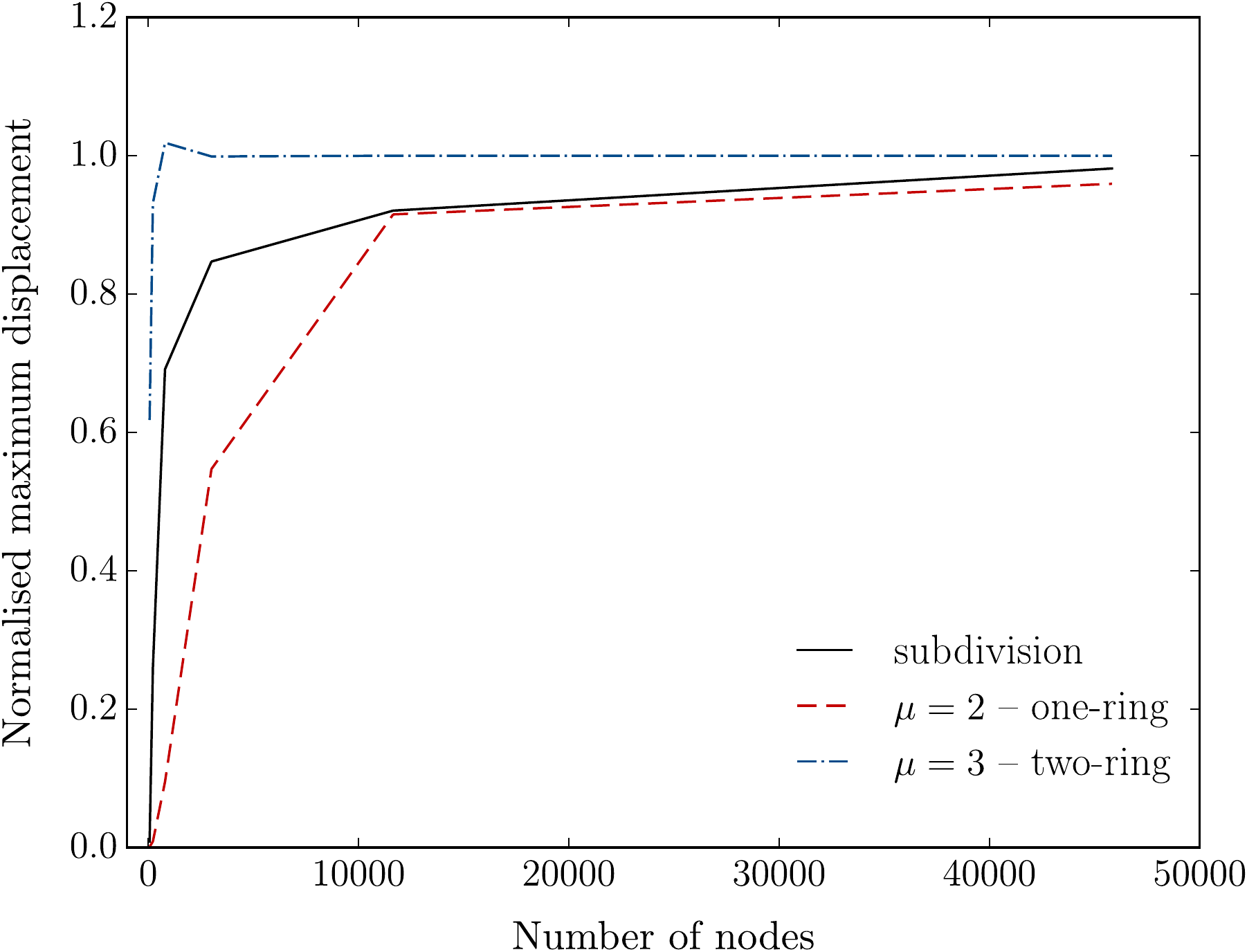}
  \caption{Pinched hemisphere. Convergence of maximum radial displacement with manifold and Catmull-Clark basis functions. The local polynomial basis~$\vec p$ is either quadratic~($\mu=2$) or cubic~($\mu=3$) and the blending functions~$w_i$  are normalised  cubic b-splines.}
  \label{fig:hemi_3}
\end{figure}

The deformed surface of the pinched hemisphere is shown in Figure
\ref{fig:hemi_1} (right). Figure \ref{fig:hemi_3} shows the convergence of the
normalised maximum radial displacement. The displacements are normalised by
$0.0924$, given in Belytschko et al.~\cite{Belytschko:1985aa}. The same plot
also includes the convergence of the finite element solution when Catmull-Clark
basis functions are used.  
As in case of pinched cylinder example, the manifold basis functions are
constructed with quadratic local polynomials~($\mu = 2$) on one-rings and cubic
polynomials~($\mu = 3$) on two-ring patches. In both cases, normalised cubic
b-splines were used as blending functions. Figure \ref{fig:hemi_3} illustrates that the
manifold basis functions constructed using cubic polynomial patches converges
slightly faster than the Catmull-Clark subdivision basis functions.

%

\section{Conclusions}

We introduced  an isogeometric analysis technique that uses manifold smooth
basis functions on quadrilateral meshes.  Manifold techniques have a long
history in computer graphics and computer aided design and several variants have
been proposed over the years.  Our implementation closely follows Ying and
Zorin~\cite{ying2004simple} and combines manifold techniques with conformal
parameterisations and the partition of unity method. The smoothness of the basis
functions is determined by the smoothness of the blending, or
partition-of-unity, functions. The approximation properties  of the basis
functions is mainly determined by the polynomial degree used in each patch. In
the presented computations, the blending function was chosen either as a
normalised linear, quadratic or cubic b-spline leading to $C^0$, $C^1$ or
$C^2$ continuous basis functions, respectively.  As patch sizes for the manifold
construction, one- or two-ring layers of elements around a vertex were
considered. The number of vertices in a patch  determines the maximum degree of
the local polynomial that can be used in the partition-of-unity interpolation.
The  polynomial coefficients in the partition of unity interpolation  are
expressed as vertex coefficients using a least-squares procedure. The finally
obtained basis functions are smooth, rational, locally supported and are
associated to  vertices in the mesh (similar to b-splines of odd degree).
The near optimal convergence of  the introduced basis  on meshes with
extraordinary vertices could be numerically confirmed.

For future research the combination of manifold basis functions  with b-splines
and the related NURBS or  subdivision surfaces  appears   especially promising. 
B-splines have several compelling properties, including refinability, that make
them ideal for geometric modelling and numerical analysis on meshes with no
extraordinary vertices.  However,  most  b-spline techniques for dealing with
extraordinary vertices, including subdivision and many $G^k$ constructions,  do
not lead to optimally convergent finite
elements~\cite{juttler2016numerical,nguyen2014comparative, nguyen2016c}. In
contrast, as numerically demonstrated manifold basis functions  yield nearly optimally
convergent finite elements independent of the connectivity of the mesh. This
suggests to use b-splines in most of the domain and to introduce manifold basis
functions only around extraordinary vertices.  In subdivision surfaces manifold
techniques have already been used to obtain $C^2$ continuity around
extraordinary vertices
\cite{levin2006modified,zorin2006constructing,antonelli2013subdivision}.
 In these three papers, instead of the conformal map  the characteristic map of
 subdivision surfaces is used to parameterise the planar patches.  The advantage
 of the characteristic map, in comparison to the conformal map, is that it
 provides a more uniform parameterisation and does not have a singularity at the
 extraordinary vertex. Additional directions for future research include  the
 mathematical convergence analysis and the proof of the linear independence of
 the introduced manifold basis function.  To this end, the large number of
 results for partition of unity methods provide a good starting point.

\appendix
\section{Implementation \label{sec:implementation}}
%
In the following we discuss the implementation of manifold basis functions focusing on data structures and algorithms. For clarity, the discussion is restricted to the case of patches~$\hat \Omega_i$ consisting of one-ring elements. This section should be read in conjunction with Section~\ref{sec:twoDconstruction} and specifically Figure~\ref{fig:mbc_8}. In line with conventional finite element implementations it is assumed that  the element matrices and vectors are assembled by iterating over elements. 
For computing the element matrices and vectors the basis function values and derivatives  at integration points are required. The basis functions depend only on the connectivity of the mesh but not on the coordinates of the vertices. Hence, they can be precomputed as part of a preprocessing step and stored for later use. 

The construction of the manifold basis functions proceeds in several steps. First, for each (non-boundary) vertex in the mesh the elements and vertices in its one-ring are identified.  Recall that we introduced one layer of ghost elements just outside the domain of interest and  that the charts which belong to the boundary vertices do not intersect the domain. The elements and vertices in an one-ring can efficiently be identified, for instance, using a half-edge data structure. 

Following the assembly of the one-rings, we endow each patch~$\hat \Omega_j$  with a blending function~$w_j (\vec \xi)$ and a local polynomial basis~$\vec p (\vec \xi)$. The degree of~$\vec p (\vec \xi)$ is chosen such that it is equal or less than the number of the vertices in the corresponding one-ring, i.e.,~$2v+1$, where~$v$ is the valence of the centre vertex.
Subsequently,  the basis functions~$ N_I(\vec \eta_Q)$ and their derivatives~$\nabla_{\vec \eta} N_I(\vec \eta_Q)$ for each quadrature (or, Gauss) point $\vec \eta_Q = (\eta^1_Q, \eta^2_Q)$ in an element~$\Box_i$  are precomputed. Each four-noded element in the mesh lies within the overlap of four patches, hence, with   each $ N_I(\vec \eta_Q)$ receiving four contributions. The contribution of one patch~$\hat \Omega_j$ to~$N_I(\vec \eta_Q)$  belonging to the vertex~$I$ is computed as follows:
\begin{enumerate}
	\item Determine the conformal coordinate~$\vec \zeta_Q$ by applying the conformal map~$z(\vec \zeta_Q)= z(\vec \eta_Q)^{4/v}$. 
	\item Identify the image of the considered element~$\Box_i$ in the patch~$\hat \Omega_j$ by comparing vertex ID's. Based on that determine the coordinate~$\vec \xi_Q$ in~$\hat \Omega_j$ by applying the rotation
	\[
		 \vec \xi_Q  = \vec R \vec \zeta_Q \,  ,
	\] 
	where the rotation~$\vec R$ is a multiple of~$2\pi / v$.
	\item Compute the contribution of patch~$\hat \Omega_j$  to basis function~$ N_I(\vec \eta_Q)$ according to~\eqref{eq:coordinatePUM2} 
	\[
		w_j(\vec \xi_Q) \sum_k^{2v+1} p_k (\vec \xi_Q) A_{kI}  \,  ,
	\]	
	where the summation is over the components in the polynomial basis~$p_k(\vec \xi_Q)$. We choose for~$p_k(\vec \xi_Q)$ either a monomial or a Lagrangian basis. The matrix~$A_{kI}$, which is the inverse of the  least-squares matrix,  does not depend on the coordinates and can be precomputed for all possible valences and stored. For computing the contribution of patch~$\hat \Omega_j$  to the derivatives~$\nabla_{\vec \eta} N_I(\vec \eta_Q)$  the above equation is differentiated, i.e., 
	\[
	\nabla_{\vec \eta} w_j = \nabla_{\vec \xi} w_j \vec R^\trans  \left ( \frac{\partial \vec \eta}{ \partial \vec \zeta} \right )^{-1}  \, .
	\] 
\end{enumerate}
After summing up the contributions of the four overlapping charts, the basis functions, their derivatives and corresponding vertices are stored in maps. The number of (non-zero) basis functions in an element is equal to the number of unique vertices in the union of the element's four charts.

%


\clearpage

\bibliographystyle{elsarticle-num-names}
\bibliography{manifoldFEM.bib}

\end{document}